\documentclass[11pt,reqno]{amsart}

\usepackage{amssymb,latexsym,verbatim,pdfsync,color,mathrsfs,graphicx,caption,epstopdf,subfigure} 

\newcommand{\bbC}{{\mathbb C}}

\newcommand{\bbR}{{\mathbb R}}
\newcommand{\bbT}{{\mathbb T}}
\newcommand{\bbZ}{{\mathbb Z}}

\def\cN{{\mathcal N}}

\def\cX{{\mathcal X}}

\def\sH{{\mathscr H}}

\def\sP{{\mathscr P}}

\def\Re{\operatorname{Re}}
\def\Im{\operatorname{Im}}

\def\la{\langle}
\def\ra{\rangle}

\def\eps{\varepsilon}
\def\vp{\varphi}
\def\ov{\overline}
\def\p{\partial}

\def\ms{\medskip}

\def\half{ {\textstyle{\frac12}}}
\def\pt{{\textstyle \frac{\pi}{2} }}

\def\Hol{\operatorname{Hol}}

\newtheorem{thm}{Theorem}[section]

\newtheorem{cor}[thm]{Corollary}
\newtheorem{lem}[thm]{Lemma}

\newtheorem{theorem}{Theorem}

\begin{document}

\title[Szeg\"o projection on worm domains]{Regularity of the Szeg\"o
  projection \\ on model worm domains}
  
\author[A. Monguzzi]{Alessandro Monguzzi}
\author[M. M. Peloso]{Marco M. Peloso}
\address{Dipartimento di Matematica ``F. Enriques''\\
Universit\`a degli Studi di Milano\\
Via C. Saldini 50\\
I-20133 Milano}
\email{alessandro.monguzzi@unimi.it}
\email{marco.peloso@unimi.it}
\thanks{Both authors supported in part by the 2010-11 PRIN grant
  \emph{Real and Complex Manifolds: Geometry, Topology and Harmonic Analysis}  
  of the Italian Ministry of Education (MIUR)}

\begin{abstract}
 In this paper we study the regularity of the  Szeg\"o projection on
 Lebesgue and Sobolev spaces on the boundary of the unbounded model
 worm domain $D'_\beta$.  
 
 We consider the Hardy space $H^2(D'_\beta)$.  Denoting by $bD'_\beta$
 the boundary of $D'_\beta$, it is classical that  $H^2(D'_\beta)$ can be
 identified with the closed subspace of $L^2(bD'_\beta,d\sigma)$,
 denoted by $H^2(bD'_\beta)$,
 consisting of the boundary values of functions in $H^2(D'_\beta)$,
 where $d\sigma$ is the induced Lebesgue measure.
 The orthogonal Hilbert space projection $P: L^2(bD'_\beta,d\sigma)\to 
 H^2(bD'_\beta)$ is called the Szeg\"o projection.
 
 Let $W^{s,p}(bD'_\beta)$ denote the Lebesgue--Sobolev space on
 $bD'_\beta$.
 We prove that $P$, initially defined on the dense subspace 
 $W^{s,p}(bD'_\beta)\cap L^2(bD'_\beta,d\sigma)$, 
  extends to a bounded operator $P: W^{s,p}(bD'_\beta)\to
  W^{s,p}(bD'_\beta)$, for $1<p<\infty$ and $s\ge0$.
\end{abstract}

\maketitle

%
\begin{center}
 {\em Dedicated to John Ryan in occasion of his 60th birthday}
 \ms
 \end{center}  

\section{Introduction and statement of the main results}
\label{sec1}

For $\beta>\pt$ we consider the domain 
\begin{equation}\label{D-beta-prime}
D'_\beta
= \Big\{ (z_1,z_2)\in\bbC^2:\, \big|\Im
z_1-\log|z_2|^2\big|<\pt,\ \big|
\log|z_2|^2\big|<\beta-\pt \Big\}
\, .
\end{equation}

Notice that $D'_\beta$ is  pseudoconvex, unbounded, and with Lipschitz
boundary.  This domain proved to be instrumental in establishing the
irregularity of the Bergman projection \cite{Barrett-Acta} and the 
failure of hypoellipticity of the $\ov\p$-Neumann problem \cite{Christ}
on the Diederich--Forn\ae ss
worm domain $\Omega_\beta$. 
The domain $\Omega_\beta$ was introduced by Diederich and Forn\ae ss 
\cite{MR0430315} as a smooth bounded pseudoconvex domain with
non-trivial {\em Nebenh\"ulle}. 
We refer the reader also to \cite{KP-JGEA} for 
an in depth discussion of properties and open problems on
$\Omega_\beta$. \ms

Among the open problems on 
$\Omega_\beta$, we now discuss the question of the (ir)-regularity of
the Szeg\"o projection and of complex Green operator on the boundary
$b\Omega_\beta$ of 
$\Omega_\beta$.  The domain $D'_\beta$ is biholomorphic to 
$$
D_\beta =
\Big\{ (w_1,w_2)\in\bbC^2:\, \Re\big(w_1e^{-i\log|w_2|^2}\big)>0,\ 
\big| \log|w_2|^2\big|<\beta-\pt \Big\} \,,
$$
via the transformation $(z_1,z_2)\mapsto (e^{z_1},z_2)$. 
D. Barrett \cite{Barrett-Acta} showed that
$\Omega_\beta$ can be exhausted by dilations of $D_\beta$ and that a
given regularity of the Bergman projection on $\Omega_\beta$ implies the 
corresponding regularity of the Bergman projection on $D_\beta$.
 He proved  the asymptotic expansion of the Bergman kernel on $D'_\beta$
and  used the transformation  rule for the Bergman kernel and
projection under biholomorphic  mappings to show  the irregularity of the
Bergman projection on $D_\beta$, thus proving the irregularity of the 
Bergman projection on $\Omega_\beta$.

Our program is to study the (ir)-regularity of the Szeg\"o projection
adapting the above paradigm. Moreover, in analogy  to M. Christ's work
\cite{Christ} and using results from \cite{HPR}, we would like to
study the hypoellipticity of the complex Green operator  associated
 to the worm domain by means of the (ir)-regularity of the Szeg\"o
projection.  There exist  several differences between the cases
of 
Bergman and Szeg\"o projections, most  noticeably the lack of a
transformation rule for the  Szeg\"o kernel 
and projection in complex dimension greater than 1.  This paper
is a step  in our program.  We provide a detailed study of the
Hardy spaces on $D'_\beta$, obtaining the expression of the Szeg\"o
kernel and proving the regularity of the Szeg\"o projection on
Lebesgue--Sobolev spaces.
\ms

We now describe our results in greater  detail.
For $\eps>0$ we consider the subdomains
$$
D'_{\beta,\eps}
= \Big\{ (z_1,z_2)\in\bbC^2:\, \big|\Im
z_1-\log|z_2|^2\big|<\pt-\eps,\ \big|
\log|z_2|^2\big|<\beta-\pt-\eps\Big\}
$$
and their 
boundaries
\begin{align*}
bD'_{\beta,\eps}
& = \Big\{ (z_1,z_2)\in\bbC^2:\, \big|\Im
z_1-\log|z_2|^2\big|=\pt-\eps,\ \big|
\log|z_2|^2\big|\leq\beta-\pt-\eps\Big\} \\
& \qquad\cup
\Big\{ (z_1,z_2)\in\bbC^2:\, \big|\Im
z_1-\log|z_2|^2\big|\leq\pt-\eps,\ \big|
\log|z_2|^2\big|=\beta-\pt-\eps\Big\} 
\,.
\end{align*} 

For $1\le p<\infty$,
we define the Hardy space of $D'_\beta$,
\begin{equation}
H^p =\Big\{ f\in\Hol(D'_\beta): \sup_{\eps>0}
\int_{bD'_{\beta,\eps}} |f|^p\, d\sigma_\eps <+\infty \Big\} \,,
\end{equation}
and set 
$$
\|f\|_{H^p}^p = \sup_{\eps>0}
\int_{bD'_{\beta,\eps}} |f|^p\, d\sigma_\eps \,.
$$

Standard basic facts of Hardy space theory give that
any $f\in H^p$ admits a boundary value function,
that we still denote by $f$, defined on 
$bD'_\beta$ that is $p$-integrable w.r.t. $d\sigma$.  Moreover, we
have the equality
$$
\|f\|_{H^p}^p 
= \int_{bD'_\beta} |f|^p\, d\sigma\,.
$$
Then, we can identify $H^p$ with a subspace of
$L^p(bD'_\beta,d\sigma)$, that is closed and that we denote by
$H^p(bD'_\beta)$.  
  The orthogonal projection
of $L^2(bD'_\beta,d\sigma)$ onto $H^2(bD'_\beta)$
$$
P: L^2(bD'_\beta,d\sigma) \to H^2(bD'_\beta)
$$
is the Szeg\"o projection  and it is described as an integral
operator, whose kernel is the Szeg\"o kernel of $D'_\beta$.

We consider the Lebesgue--Sobolev spaces
 $W^{s,p}(bD'_\beta)$ on $bD'_\beta$, defined for $s\ge0$ and $1\le
p\le \infty$.  We observe that $bD'_\beta$ is a Lipschitz embedded
hypersurface in $\bbC^2$ and therefore these spaces are defined
by the standard well-known theory.  More simply, $bD'_\beta$ is union
of 4 disjoint connected components and a set of surface measure zero,
the {\em distinguished boundary},
and each of these components of positive measure
 is in fact a portion of a hyperplane in
$\bbR^4$. 

Our main results are the following.  
\begin{theorem}\label{Sze-ker-main-thm}
The Szeg\"o kernel of $D'_\beta$ is 
given by
\begin{equation}\label{Sze-ker}
K(z,w) = \sum_{j\in\bbZ}  {z_2}^j{\ov w_2}^j k_j(z_1,w_1) \,,
\end{equation}
where $z=(z_1,z_2)$, $w=(w_1,w_2)\in D'_\beta$ and
$$
k_j (z_1,w_1) = \frac{1}{2\pi} 
 \int_\bbR  \frac{e^{i(z_1-\ov w_1)\xi} }{\nu(\xi,j)}\, d\xi\,,
$$
\begin{multline*}
\nu(\xi,j) 
= (\beta-\pt) \cosh (\pi\xi)
\int_{-1}^1 e^{-(\beta-\frac\pi2)(2\xi-(j+1))s}
\sqrt{1+4e^{-(\beta-\frac\pi2)s}}\, ds  \\+  
2\cosh \big[  (\beta-\pt)(2\xi-(j+\half)) \big]
\frac{\sinh (\pi\xi)}{\xi} \,.
\end{multline*}

The series in \eqref{Sze-ker} converges in $H^2(D'_\beta)$ for every
$(w_1,w_2)$ fixed,
and
uniformly in compact subsets of $D'_\beta\times D'_\beta$. 
\end{theorem}

\begin{theorem}\label{Sze-pro-main-thm}
The Szeg\"o projection $P$, initially defined
on the dense subspace 
$W^{s,p}(bD'_\beta)\cap L^2(bD'_\beta,d\sigma)$, 
 extends to a bounded operator 
$$
P: W^{s,p}(bD'_\beta)\to
 W^{s,p}(bD'_\beta)\,, 
$$
for $1<p<\infty$ and $s\ge0$.
\end{theorem}
\ms

Our results thus  concern the  Lebesgue--Sobolev regularity of the Szeg\"o projection
on certain domains in $\bbC^2$.  Of course this is a classical and
widely studied problem and there exists a vast literature.  
Hence, we restrict this discussion to the case of domains in $\bbC^n$,
with $n>1$.   The cases  of smoothly bounded strongly pseudoconvex domains, 
domains of finite type in $\bbC^2$, and convex domains of finite type
in  in $\bbC^n$ are classical by now, and we refer to
\cite{Phong-Stein}, \cite{NRSW-Annals},  and \cite{McN-St},  
 resp. 
On all such domains $\Omega$, the Szeg\"o projection is bounded
on the Lebesgue $W^{s,p}(b\Omega)$ for $1<p<\infty$  and $s\ge0$, and it
is also weak-type $(1,1)$.  These results were extended by K. Koenig
\cite{MR1879002} to domains in CR manifolds on which the tangential
Cauchy--Riemann operator has closed range and satisfying some conditions
on the eigenvalues of the Levi form.

More recently,
L. Lanzani and E. M. Stein studied the regularity of
the Szeg\"o projection on strongly pseudoconvex domains with minimal
boundary regularity, \cite{LS1,LS2}, obtaining the $L^p$-boundedness
for $1<p<\infty$.  \ms

The domain $D'_\beta$ is unbounded and its boundary consists of
hyperplanes, hence it is Levi flat.  On such domain, S. Krantz and the
second author proved that the Bergman projection is bounded on $L^p$
for $1<p<\infty$ \cite{KP-Houston} and the first author extended to
the scale of Sobolev spaces \cite{Monguzzi2015}.
As mentioned above, $D'_\beta$ is biholomorphically equivalent to another
model worm domain $D_\beta$.   On $D_\beta$, in contrast,  
we expect the Szeg\"o projection to show some irregularity, as in the
case of the Bergman  projection. In principle this should indicate that 
the Szeg\"o projection on the smooth worm domain $\Omega_\beta$
would not preserve Sobolev spaces $W^s$ for $s\geq s_0$ where $s_0$ is
a positive real number depending on the geometry of $\Omega_\beta$,
namely, $s_0\to 0$ as $\beta\to+\infty$.  These 
investigations are part of our program and are going to be topics of future
work.  
\ms

We mention that on the domains $D'_\beta$ and $D_\beta$ it is possible to
introduce and study also Hardy spaces defined in terms of the induced 
surface measure on the {\em distinguished boundary}, 
that is, in the case of the domain $D'_\beta$, the four
vertices in Figure \ref{figure1}.  This has been done by the first author
\cite{MR3440104} in the case of $D'_\beta$, and in \cite{MP-in-prep}
in the case of $D_\beta$.   More precisely, denote by $\p D'_\beta $
and $\p D_\beta$ the distinguished boundaries of $D'_\beta $
and $ D_\beta$, resp., and by $\sP'$ and $\sP$ the corresponding
projections, resp.  We point out that in this setting, the operators
$\sP'$ and $\sP$ are given by singular
integrals over  $\p D'_\beta $
and $\p D_\beta$, resp.
Then, 
the following results are proved: $\sP': W^{s,p} (\p D'_\beta) \to
W^{s,p} (\p D'_\beta) $, when $1<p<\infty$ and $s\ge0$ (\cite{MR3440104}), while 
$\sP: L^p (\p D_\beta) \to
L^p  (\p D_\beta) $, if and only if $p\in \big(\frac{2}{1+\nu_\beta},
\frac{2}{1-\nu_\beta} \big)$, where $\nu_\beta=\pi/(2\beta-\pi)$ (\cite{MP-in-prep}). 
\ms

It would also be interesting studying the regularity of the Szeg\"o
projection on domains that are models for the higher-dimensional
worm domains introduced by Barrett and {\c{S}}ahuto{\u{g}}lu in
\cite{B-S}.  

 For other recent papers related to the
topic of this work we also mention \cite{KPS, BEP, KPS2}.
\ms 

The plan of the paper is the following.  In Section \ref{2} we obtain
a representation of the Szeg\"o kernel of $D'_\beta$ and prove Theorem
\ref{Sze-ker-main-thm}.  In Section \ref{Lp-bdb-sec} we decompose the
Szeg\"o projection as sum of operators $T_{\ell,\ell'}$,
$\ell,\ell'\in\{1,\dots,4\}$,  and prove their
boundedness on $L^p$ spaces. We find particularly interesting that the 
operators $T_{\ell,\ell'}$ can be written as composition of Fourier
multiplier operators on $\bbR\times\bbT$ and an integral operator of
Hilbert type.  We can then apply a combination of the classical Marcinkiewicz
multiplier theorem  and Schur's test  
to obtain the $L^p$-boundedness of $P$.  
In Section \ref{Lemmas-Sec} we complete the proof of Theorem
\ref{Sze-pro-main-thm} and prove a few remaining technical lemmas.  
\ms

\noindent
{\em Acknowledgments.}
We laid the groundwork of this project when we were both visiting the
Department of Mathematical Sciences of University of Arkansas.  We
wish to thank such institution for the great hospitality
and for providing a very stimulating working environment.

The second author wishes to thank in particular John Ryan for the
many interesting and pleasant meetings and discussions had during
his stay at the Department of Mathematical Sciences.

We also thanks the anonymous reviewers for carefully reading the
manuscripts and making a number of remarks that helped to improve the
clarity of the paper. 

\ms

\section{The Hardy space and the Szeg\"o kernel}\label{2}

\subsection{Decomposition of the Hardy space}

In order to describe the norm on the Hardy spaces, we begin by
computing 
$d\sigma_\eps$.   We label the sides of $bD'_{\beta,\eps}$ as
 $E_{1,\eps},\dots,E_{4,\eps}$, where $E_{1,\eps}$ is
the right side, and then proceeding counter-clockwise.  Explicitly:
\begin{align*}
E_{1,\eps} 
& =  \big\{ (z_1,z_2)\in\bbC^2:\, \Im
z_1-\log|z_2|^2 =\pt-\eps,\ \big|
\log|z_2|^2\big|\le\beta-\pt-\eps\big\}\,; \\
E_{2,\eps} 
& = \big\{ (z_1,z_2)\in\bbC^2:\, \big|\Im
z_1-\log|z_2|^2\big|\le\pt-\eps,\  
\log|z_2|^2 =\beta-\pt-\eps\big\}\,; \\
E_{3,\eps}
& =  \big\{ (z_1,z_2)\in\bbC^2:\, \Im
z_1-\log|z_2|^2 =-\pt+\eps,\ \big|
\log|z_2|^2\big|\le\beta-\pt-\eps\big\}\,; \\
E_{4,\eps}
& = \big\{ (z_1,z_2)\in\bbC^2:\, \big|\Im
z_1-\log|z_2|^2\big|\le\pt-\eps,\  
\log|z_2|^2 =-(\beta-\pt-\eps)\big\}\,.
\end{align*}

\begin{figure}[!h]
\begin{center}
\includegraphics[width=15cm]{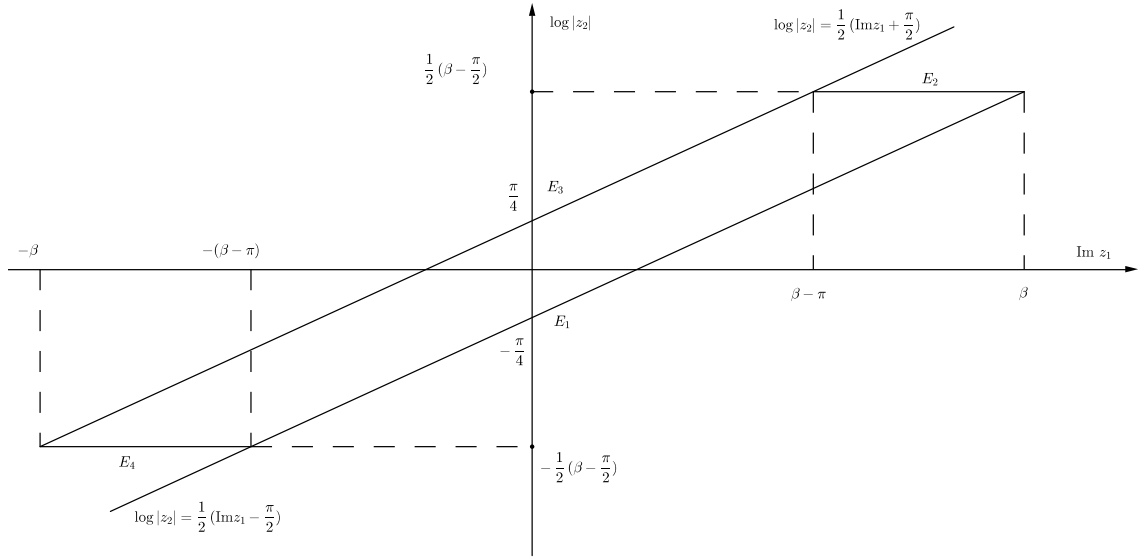}
\end{center}
\caption{A representation of the domain $D'_\beta$ in the $(\Im z_1,\log|z_2|)$-plane.}\label{figure1}
\end{figure}

We also use the analogous decomposition on $bD'_\beta$ and denote the
four sides by $E_1,\dots,E_4$. See Figure \ref{figure1}.

We wish to determine the induced surface measure on
$bD'_{\beta,\eps}$, for $\eps\ge0$.  It is well known that, if $M$ is
an embedded manifold in $\bbR^n$, parametrized by $F:A\subset\bbR^m\to
\bbR^n$, then 
$$
\int_M f d\sigma = \int_A f(F(s)) \sqrt{\operatorname{det} ({\,}^t(JF) JF) (s)
}\, ds \,,
$$
where $JF$ denotes the jacobian matrix of $F$.

In the case of $E_{1,\eps}$ and $E_{3,\eps}$, resp., we have 
$$
F(x,r,\theta) = \begin{pmatrix}
x \\ \log r^2+\pt-\eps \\  r\cos\theta \\ r\sin\theta
\end{pmatrix},
\quad
\text{and }\ F(x,r,\theta) = \begin{pmatrix}
x \\ \log r^2-\pt+\eps \\  r\cos\theta \\ r\sin\theta
\end{pmatrix} , \text{resp.}
$$
 with
domain $A=\bbR\times\big(e^{-\frac12(\beta-\frac\pi2-\eps)},
e^{\frac12(\beta-\frac\pi2-\eps)} \big)\times(0,2\pi)$, 
so that  $d\sigma =dx\, r\sqrt{1+4r^{-2}}dr\, d\theta$
(in both cases).

Analogously, in the case of
$E_{2,\eps}$ and $E_{4,\eps}$, resp., we have
$$
F(x,y,\theta) = \begin{pmatrix}
x \\ y  \\  e^{\frac12(\beta-\frac\pi2-\eps)} \cos\theta \\
e^{\frac12(\beta-\frac\pi2-\eps)} \sin\theta
\end{pmatrix}
\quad
\text{and }\ F(x,y,\theta) = \begin{pmatrix}
x \\ y\\  e^{-\frac12(\beta-\frac\pi2-\eps)}\cos\theta \\ e^{-\frac12(\beta-\frac\pi2-\eps)}\sin\theta
\end{pmatrix} , \text{resp.}
$$
 with
domain $A=\bbR\times\big\{ |y-(\beta-\frac\pi2-\eps)|<\frac\pi2 \big\}
\times(0,2\pi)$, and 
$A=\bbR\times\big\{ |y+(\beta-\frac\pi2-\eps)|<\frac\pi2 \big\}
\times(0,2\pi)$, resp. 
Therefore,  $d\sigma =dx\, dy\, e^{\frac12(\beta-\frac\pi2-\eps)}\, d\theta$
and $d\sigma =dx\, dy\, e^{-\frac12(\beta-\frac\pi2-\eps)}\, d\theta$,
resp.

We observe that the four vertices of the parallelogram in Figure
\ref{figure1} have zero surface measure and therefore can be
disregarded and we shall do so in what follows.    

\ms

Therefore, $f\in H^p$ if and only if $\| f\|_{H^p}^p
 = \sup_{\eps>0} \int_{bD'_{\beta,\eps}}|f|^pd\sigma_\eps<\infty$,
 where
\begin{align*}
\int_{bD'_{\beta,\eps}}|f|^pd\sigma_\eps
& = \bigg( \int_{E_{1,\eps}}  +\int_{E_{2,\eps}} + \int_{E_{3,,\eps}}  +
\int_{E_{4,\eps}} \bigg) |f|^pd\sigma_\eps \\
& = \int_0^{2\pi}
\int_{e^{-\frac12(\beta-\frac\pi2-\eps)}}^{e^{\frac12(\beta-\frac\pi2-\eps)}}
\int_\bbR |f\big(x+i(\log r^2+\pt-\eps),re^{i\theta}\big)|^p\,dx
\sqrt{1+4r^{-2}}\, rdr\, d\theta \\
& \quad + \int_0^{2\pi}
\int_\bbR
\int_{ |y-(\beta-\frac\pi2-\eps)|<\frac\pi2} 
 |f\big(x+iy, e^{\frac12(\beta-\frac\pi2-\eps)}e^{i\theta}
 \big)|^p\,dy dx\, e^{\frac12(\beta-\frac\pi2-\eps)}\, d\theta  \\
& \quad \quad + \int_0^{2\pi}
\int_{e^{-\frac12(\beta-\frac\pi2-\eps)}}^{e^{\frac12(\beta-\frac\pi2-\eps)}}
\int_\bbR |f\big(x+i(\log r^2-\pt+\eps),re^{i\theta}\big)|^p\,dx
\sqrt{1+4r^{-2}}\, rdr\, d\theta \\ 
& \quad \quad \quad + \int_0^{2\pi}
\int_\bbR
\int_{ |y+(\beta-\frac\pi2-\eps)|<\frac\pi2} 
 |f\big(x+iy, e^{-\frac12(\beta-\frac\pi2-\eps)}e^{i\theta}
 \big)|^p\,dy dx\, e^{-\frac12(\beta-\frac\pi2-\eps)} d\theta \,. \ms
\end{align*}

Using standard one-variable results, it is easy to see that the above
sup is attained as $\eps\to0^+$.  Denoting by $f$ also the boundary
value function of $f$, we then have
\begin{align} 
\|f\|_{H^p}^p
& = \int_0^{2\pi}
\int_{e^{-\frac12(\beta-\frac\pi2)}}^{e^{\frac12(\beta-\frac\pi2)}}
\int_\bbR |f\big(x+i(\log r^2+\pt),re^{i\theta}\big)|^p\,dx\,
r\sqrt{1+4r^{-2}} \,dr d\theta  \notag \\
& \quad + e^{\frac12(\beta-\frac\pi2)}  \int_0^{2\pi}
\int_\bbR
\int_{ |y-(\beta-\frac\pi2)|<\frac\pi2} 
 |f\big(x+iy, e^{\frac12(\beta-\frac\pi2)}e^{i\theta}
 \big)|^p\,dy dx d\theta\notag \\
& \quad \quad + \int_0^{2\pi}
\int_{e^{-\frac12(\beta-\frac\pi2)}}^{e^{\frac12(\beta-\frac\pi2)}}
\int_\bbR |f\big(x+i(\log r^2-\pt),re^{i\theta}\big)|^p\,dx\,
r\sqrt{1+4r^{-2}}\, dr\, d\theta  \notag \\
& \quad \quad \quad +e^{-\frac12(\beta-\frac\pi2)}  \int_0^{2\pi}
\int_\bbR
\int_{ |y+(\beta-\frac\pi2)|<\frac\pi2} 
 |f\big(x+iy, e^{-\frac12(\beta-\frac\pi2)}e^{i\theta}
 \big)|^p\,dy dx d\theta \notag \\
& =  \half \int_0^{2\pi}
\int_{-\beta+\pi}^{\beta} 
\int_\bbR |f\big(x+iy ,e^{\frac12(y-\frac\pi2)}e^{i\theta}\big)|^p\,dx
\, e^{y-\frac\pi2} \sqrt{1+4e^{-(y-\frac\pi2)}}\, dyd\theta  \notag \\
& \quad + e^{\frac12(\beta-\frac\pi2)} \int_0^{2\pi}
 \int_{ \beta-\pi}^\beta \int_\bbR 
 |f\big(x+iy, e^{\frac12(\beta-\frac\pi2)}e^{i\theta}
 \big)|^p\, dxdyd\theta   \notag \\
& \quad \quad + \half \int_0^{2\pi}
\int_{-\beta}^{\beta-\pi} 
\int_\bbR |f\big(x+iy ,e^{\frac12(y+\frac\pi2)} e^{i\theta}\big)|^p\,dx\,
e^{y+\frac\pi2} \sqrt{1+4e^{-(y+\frac\pi2)}}\, dy d\theta \notag \\
& \quad \quad \quad + e^{-\frac12(\beta-\frac\pi2)} \int_0^{2\pi}
 \int_{ -\beta}^{-\beta+\pi} \int_\bbR 
 |f\big(x+iy, e^{-\frac12(\beta-\frac\pi2)}e^{i\theta}
 \big)|^p\,dxdyd\theta \,. \ms \label{norm1}
\end{align}

We consider the Hilbert space $H^2$ and  decompose it  following
\cite{Ki}, see also \cite{Barrett-Acta,KP-Houston}.
Using Fourier series
expansion in the second variable we have the orthogonal decomposition  
$$
H^2 = \bigoplus_{j\in\bbZ} \sH_j
$$
where for $j\in\bbZ$, 
$\sH_j$ is the subspace of functions $f$ in $H^2$ such that
$f(z_1,e^{i\theta}z_2)= e^{ij\theta} f(z_1,z_2)$.  The orthogonal
projection $Q_j: H^2\to \sH_j$ is  
$$
Q_j f(z_1,z_2)=
\frac{1}{2\pi} \int_0^{2\pi}
f(z_1,z_2e^{i\theta}) e^{-ij\theta}\, d\theta \,.
$$
The function
$$
\bigg(\frac{1}{2\pi } \int_0^{2\pi}
f(z_1,z_2e^{i\theta}) e^{-ij\theta}\, d\theta \bigg)  z_2^{-j} 
$$
 
is holomorphic  in $D'_\beta$ and locally constant in $|z_2|$.  Since the
fibers over $z_1$ are connected, it is actually constant in $z_2$.  Hence,  we
simply write $f_j(z_1)$, and notice that $f_j\in \Hol(S_\beta)$, 
where 
$$
S_\beta=\big\{x+iy\in\bbC: |y|<\beta\big\}\,.
$$
  
Thus, 
$$
\sH_j = \Big\{ f\in H^2:\, f(z_1,z_2)=f_j(z_1)z_2^j,\ \text{\ with\ }
f_j\in \Hol(S_\beta) \Big\}\,.
$$

Now, let $f(z_1,z_2)=\sum_{j\in\bbZ} f_j(z_1)z_2^j$. By orthogonality
we have 
$
\| f\|_{H^2}^2 
= \sum_{j\in\bbZ}  \| f_j z_2^j\|_{H^2}^2 $.  From \eqref{norm1} we compute
\begin{align*}
\| f_j z_2^j\|_{H^2}^2 
& =  \pi \int_{-\beta+\pi}^{\beta} 
\int_\bbR |f_j(x+iy)|^2\,dx
\, e^{(y-\frac\pi2)(j+1)} \sqrt{1+4e^{-(y-\frac\pi2)}} \,dy \notag \\
& \quad + 
2\pi e^{(\beta-\frac\pi2)(j+\frac12)}  
\int_\bbR 
\int_{\beta-\pi}^\beta 
 |f_j (x+iy)|^2\,dy dx 
 \notag \\
& \quad \quad 
+ \pi
\int_{-\beta}^{\beta-\pi} 
\int_\bbR |f_j(x+iy)|^2\,dx \,
e^{(y+\frac\pi2)(j+1)}\sqrt{1+4e^{-(y+\frac\pi2)}}\, dy
\notag \\
& \quad \quad \quad+ 2\pi e^{-(\beta-\frac\pi2)(j+\frac12)} 
\int_\bbR 
\int_{-\beta}^{-\beta+\pi} 
  |f_j (x+iy)|^2\,dy dx \\
& = \int_{S_\beta} |f_j(z_1)|^2\, \omega_j(y) dA(z_1) \,,
\end{align*}
where $z_1=x+iy$, and
\begin{multline}\label{omega-j}
\omega_j (y)
=  \pi \chi_{(-\beta+\pi,\beta)}(y)  e^{(y-\frac\pi2)(j+1)} \sqrt{1+4e^{-(y-\frac\pi2)}}
+ 2 \pi \chi_{(\beta-\pi,\beta)}(y) 
e^{(\beta-\frac\pi2)(j+\frac12)} \\
  + \pi \chi_{(-\beta,\beta-\pi)}(y)  e^{(y+\frac\pi2)(j+1)} \sqrt{1+4e^{-(y+\frac\pi2)}} + 
 2\pi \chi_{(-\beta,-\beta+\pi)}(y) 
e^{-(\beta-\frac\pi2)(j+\frac12)}  \,. 
\end{multline}

Thus, $f_j$ is not merely holomorphic in $S_\beta$, it belongs to
weighted Bergman space
$A_j^2=A^2(S_\beta, \omega_j(y))dA$.  Since $\omega_j$ is bounded
above and below by positive constant, 
$A_j^2$ 
coincides with the unweighted Bergman space of $S_\beta$.  \ms

In order to obtain the expression of Szeg\"o kernel we wish to
compute the reproducing kernel for $A^2_j$.  \ms

\subsection{The weighted Bergman space of the strip}
\ms

We first prove a Paley--Wiener type theorem for $A^2_j$. 
\begin{thm}\label{PW-thm} {\rm {\bf (1)}}
Let $\omega_j$ be as in \eqref{omega-j}.
Then
\begin{multline}\label{hat-omega-j}
\widehat \omega_j(-2i\xi)= 2\pi (\beta-\pt) 
\cosh (\pi\xi)
\int_{-1}^1 e^{-(\beta-\frac\pi2)(2\xi-(j+1))s}
\sqrt{1+4e^{-(\beta-\frac\pi2)s}}\, ds  \\+  4\pi
\cosh \big[  (\beta-\pt)(2\xi-(j+\half) \big]
\frac{\sinh (\pi\xi)}{\xi} 
 \,.
\end{multline}

{\rm {\bf (2)}} Let $f\in A^2_j$.  Then $f_y=f(\cdot+iy)
\in L^2(\bbR)$ for a.e. $y\in(-\beta,\beta)$, $\widehat{f_y}  (\xi)=
e^{-y\xi}\widehat f_0 (\xi)$, $\widehat f_0\in L^2(\bbR, \nu(\xi,j)d\xi)$ and
\begin{equation}\label{Plan-S-beta}
\|f\|_{A_j^2}^2 
= \int_\bbR |\widehat f_0(\xi)|^2\, \nu(\xi,j)d\xi\,,
\end{equation}
where $\nu(\xi,j) =\frac{1}{2\pi} \widehat\omega_j(-2i\xi)$.

{\rm {\bf (3)}} 
Conversely, let $\omega_j$ be as in \eqref{omega-j} and let 
$\nu(\xi,j) =\frac{1}{2\pi} \widehat\omega_j(-2i\xi)$.  Then, given
$g\in  L^2(\bbR, \nu(\xi,j)d\xi)$  for $z\in S_\beta$ we define
\begin{equation}\label{conv-PW}
f(z) =\frac{1}{2\pi} \int_\bbR e^{iz\xi} g(\xi)\, d\xi\,.
\end{equation}
Then $f\in A^2_j$, $\widehat{f_0}=g$ and \eqref{Plan-S-beta}
holds. 
\end{thm}

\proof
{\bf (1)} This is an elementary calculation.  Decomposing $\omega_j$ as in 
\eqref{omega-j} as sum of four terms we set
$$
\frac{1}{2\pi}  \widehat \omega_j(-2i\xi)   
= \frac{1}{2\pi}  \int_{-\beta}^\beta  e^{-2y\xi} \omega_j(y)\, dy 
= I_1+\dots+I_4\,.
$$
Then,
\begin{align*}
I_1+I_3
& = \half\! \int_{-\beta+\pi }^{\beta}  e^{-2y\xi} 
  e^{(y-\frac\pi2)(j+1)} \sqrt{1+4e^{-(y-\frac\pi2)}} \,dy 
+ \half\! \int_{-\beta}^{\beta-\pi}    e^{-2y\xi} 
e^{(y+\frac\pi2)(j+1)} \sqrt{1+4e^{-(y+\frac\pi2)}} \,dy \\
& = \cosh (\pi\xi)
\int_{-\beta+\frac\pi2}^{\beta-\frac\pi2}  e^{-2t\xi} e^{t(j+1)}
\sqrt{1+4e^{-t}}\, dt \\
& = \cosh (\pi\xi)
\int_{-1}^1 e^{-(\beta-\frac\pi2)(2\xi-(j+1))s}
\sqrt{1+4e^{-(\beta-\frac\pi2)s}}\, ds 
 \,,
\end{align*}
and
\begin{align*}
I_2+I_4
& = e^{(\beta-\frac\pi2)(j+\frac12)} 
\int_{\beta-\pi}^\beta e^{-2y\xi}  \, dy 
+  e^{-(\beta-\frac\pi2)(j+\frac12)} 
\int_{-\beta}^{-\beta+\pi} e^{-2y\xi}  \, dy \\
& =  e^{-(\beta-\frac\pi2)(2\xi-(j+\frac12)} \frac{\sinh (\pi\xi)}{\xi} 
+  e^{(\beta-\frac\pi2)(2\xi-(j+\frac12))} \frac{\sinh (\pi\xi)}{\xi} \\
& = 2\cosh \big[  (\beta-\pt)(2\xi-(j+\half)) \big]
\frac{\sinh (\pi\xi)}{\xi} 
 \,.
\end{align*}
This proves {\bf (1)}.
\ms

{\bf (2)} We have already observed
that $f\in A_j^2$ if and only if $f\in
A^2(S_\beta)$.  Then $f_y=f(\cdot+iy)
\in L^2(\bbR)$ for a.e. $y\in(-\beta,\beta)$, $\widehat{f_y}  (\xi)=
e^{-y\xi}\widehat f_0 (\xi)$, see 
\cite{Harper}.
Therefore, 
\begin{align*}
\| f\|_{A^2_j}^2 
& = \frac{1}{2\pi} \int_{-\beta}^\beta \int_\bbR |\widehat{f_0} (\xi)|^2 e^{-2y\xi}
\, d\xi \, \omega(y)\, dy\\
& = \frac{1}{2\pi} \int_\bbR  |\widehat{f_0}(\xi)|^2  \int_{-\beta}^\beta
e^{-2y\xi} \, \omega(y)\, dy \, d\xi \\
& = \frac{1}{2\pi} \int_\bbR  |\widehat{f_0} (\xi)|^2  \int_{-\beta}^\beta
e^{-2y\xi} \, \omega(y)\, dy \, d\xi = \int_\bbR  |\widehat{f_0} (\xi)|^2 \,
\nu(\xi,j)\, d\xi\,.
\end{align*}

{\bf (3)}
Let $g$ be continuous with compact support in $\bbR$, so that 
$g\in  L^2(\bbR, \nu(\xi,j)d\xi)\cap
L^2(\bbR, \frac{\sinh(2\beta\xi)}{\xi}d\xi)$.
Then, by \cite{Harper}, if
$f$ is defined as in \eqref{conv-PW} then $f$ is in $A^2(S_\beta)$,
 hence in $A^2_j$.  Thus, part {\bf (2)} applies and therefore
$\widehat{f_0}=g$ and \eqref{Plan-S-beta}
holds.  Since this holds true on a dense subset of 
$L^2(\bbR, \nu(\xi,j)d\xi)$ the conclusion follows.
\qed
\ms

\begin{cor}\label{rep-ker-A2j}
The reproducing kernel of $A^2_j$ is given by
\begin{equation}
k_j(z_1,w_1) = \frac{1}{2\pi} \int_\bbR 
\frac{e^{i(z_1-\ov w_1)\xi} }{\nu(\xi,j)}\, d\xi\,, \label{k-omega-j}
\end{equation}
where 
\begin{align}
\nu(\xi,j) 
& = \frac{1}{2\pi} \widehat \omega_j(-2i\xi)  \,, \label{def-nu}
\end{align}
and $\widehat \omega_j(-2i\xi)$ is given by \eqref{hat-omega-j}.
\end{cor}

Notice that $\nu(\cdot,j)$ is even, positive and $1/\nu(\cdot,j)$ is a Schwartz function
such that $1/\nu(\xi,j)\le c e^{-2\beta|\xi|}$ for a positive constant
$c$, independent of $j$, so that
$k_j$ is well defined for $z_1,w_1\in S_\beta$.
\ms

\proof Writing $k_{w_1}= k_j(\cdot,w_1)$,
where $w_1 \in S_\beta$, 
for 
$f\in A_j^2$ we
have 
$$
f(w_1) = \la f,k_{w_1}\ra_{A_j^2} = \la \widehat{f_0},\, \widehat{k_{w_1,0}}
\ra_{L^2(\nu(\xi,j)d\xi)} \,,
$$
and also 
$$
f(w_1) =\frac{1}{2\pi} \int_\bbR e^{iw_1\xi} \widehat{f_0}(\xi) \,
d\xi\,. 
$$
Therefore, $\widehat{k_{w_1,0}}(\xi)= \nu(\xi,j)^{-1} e^{\ov w_1\xi}$, that
is, \eqref{k-omega-j} holds.
\qed

\ms

\proof[Proof of Theorem \ref{Sze-ker-main-thm}]
This follows at once from the orthogonal decomposition $H^2 =
\bigoplus_{j\in\bbZ} \sH_j$, where each $\sH_j$ is isometrically
equivalent to $A_j^2$ via the correspondence $\sH_j 
\ni f_j(z_1)z_2^j \mapsto f_j\in A^2_j$, and Corollary \ref{rep-ker-A2j}.
\qed

\ms

\section{$L^p$-boundedness of the Szeg\"o projection}
\label{Lp-bdb-sec}

The Szeg\"o projection is the Hilbert space orthogonal projection of 
$L^2(bD'_\beta)$ onto the closed subspace consisting of the boundary
values of functions in $H^2$.  By standard Hilbert space theory, the 
Szeg\"o projection of 
 $\vp\in L^2(bD'_\beta)$, 
is given by 
\begin{equation*}
P\vp (\kappa)
= \int_{bD'_\beta}\vp(\kappa')  
K(\kappa,\kappa')\, d\sigma(\kappa') \,,
\end{equation*}
with $\kappa \in bD'_\beta$. 
\ms

We decompose $P$ writing
\begin{align*}
P 
& = (\chi_{E_1}+\cdots+\chi_{E_4}) P (\chi_{E_1}+\cdots+\chi_{E_4}) 
= \sum_{\ell,\ell'=1}^4 \chi_{E_{\ell}} P \chi_{E_{\ell'}}
\end{align*}
and set, for $\ell,\ell'=1,\dots,4$,
\begin{equation}\label{P-ell-ell-prime}
P_{\ell,\ell'} =\chi_{E_{\ell}} P \chi_{E_{\ell'}} \,.
\end{equation}

Thus, it suffices to study the boundedness of the
operators  $P_{\ell,\ell'}$, for $\ell,\ell'=1,\dots,4$.
We introduce the (global) coordinates on $E_1,\dots,E_4$ in order to
write down the operators $P_{\ell,\ell'}$ explicitly. 
We set
\begin{equation}\label{I-ell}
 I_1=  (-\beta+\pi,\beta),\quad I_2 =(\beta-\pi,\beta) ,\quad 
I_3 =(-\beta,\beta-\pi),\quad I_4 = (-\beta,-\beta+\pi)
  \,, 
\end{equation}
and
\begin{equation}\label{cal-X-ell}
\cX_\ell=\bbR\times  I_\ell\times\bbT,\quad \ell=1,\dots,4.
\end{equation}
Here and in what follows, we denote by $\bbT$ the 1-dimensional torus
$\bbR/2\pi\bbZ$.
\ms

The $L^p$-boundedness of the Szeg\"o projection on
$D'_\beta$ is a consequence of the following theorem.

\begin{thm}\label{Lp-bddness}
For all $\ell,\ell'\in\{1,\dots,4\}$ the following holds true.

{\rm {\bf (1)}} There exists an operator $T_{\ell,\ell'}$ such that, for $1<p<\infty$,
$P_{\ell,\ell'}: L^p(E_{\ell'})\to L^p(E_\ell)$ is bounded if and only 
$T_{\ell,\ell'}: L^p(\cX_{\ell'})\to  L^p(\cX_{\ell})$ is bounded. 

{\rm {\bf (2)}} There exists a 
Fourier multiplier operator
$M^{\ell,\ell'}_{(y,t)}$ on $\bbR\times\bbT$ depending on the
parameters $t\in I_{\ell'}$ and $y\in I_\ell$ such that 
$$
T_{\ell,\ell'} = \int_{I_{\ell'}} M^{\ell,\ell'}_{(y,t)} \,dt\,.
$$

{\rm {\bf (3)}} The operator 
$$
T_{\ell,\ell'}: L^p(\cX_{\ell'})\to  L^p(\cX_{\ell})
$$
 is bounded for $1<p<\infty$.
\end{thm}

\proof
{\bf (1)} 
This step simply consists of writing the 
operators $P_{\ell,\ell'}$ in coordinates and  simplifying  the innocuous
weight factor. \ms

For $\ell=1,\dots,4$ and $\vp \in L^p(E_\ell)$ we define
\begin{align}
\begin{split}
& \Lambda_1\vp(x,y,\theta) =
\frac{1}{\sqrt{1+4e^{-(y-\frac\pi2)}}}\vp\big(x+iy,e^{\frac12(y-\frac\pi2)}e^{i\theta}\big)
\,, \ \text{for\ }(x,y,\theta)\in\cX_1\,;\qquad
\\
& \Lambda_2\vp(x,y,\theta) =
e^{-\frac12(\beta-\frac\pi2)}\vp\big(x+iy,e^{\frac12(\beta-\frac\pi2)} e^{i\theta}\big)
\,, \ \text{for\ } (x,y,\theta)\in\cX_2\,;\qquad
\\
& \Lambda_3\vp(x,y,\theta) =
\frac{1}{\sqrt{1+4e^{-(y+\frac\pi2)}}}\vp\big(x+iy,e^{\frac12(y+\frac\pi2)}e^{i\theta}\big)
\,, \ \text{for\ }(x,y,\theta)\in\cX_3\,;\qquad 
\\
& \Lambda_4\vp(x,y,\theta) =
e^{\frac12(\beta-\frac12)}\vp\big(x+iy,e^{-\frac12(\beta-\frac\pi2)}e^{i\theta}\big)
\,, \ \text{for\ } (x,y,\theta)\in\cX_4\,.
\end{split}
\label{Lambda-ell}
\end{align}

Notice that $\Lambda_\ell:L^p(E_\ell)\to L^p(\cX_\ell)$ is an onto 
isomorphism.  Hence, $ P_{\ell,\ell'}: L^p(E_\ell)\to L^p(E_{\ell'})$ is bounded if and
only if $\Lambda_\ell P_{\ell,\ell'}\Lambda_{\ell'}^{-1} :L^p(\cX_{\ell'})\to L^p(\cX_\ell)$ is
bounded.

We set
\begin{equation}\label{T-ell-ell-prime-def}
T_{\ell,\ell'} = \Lambda_\ell P_{\ell,\ell'}\Lambda_{\ell'}^{-1} \,.
\end{equation}

{\bf (2)} We begin with $T_{1,1}$.  We denote by $C_c(\cX_\ell)$ the
continuous functions with compact support in $\cX_\ell$.  
Observe that, for $\psi \in C_c(\cX_1)$, 
$T_{1,1}$ is given
by 
\begin{align*}
T_{1,1}\psi (x,y,\theta)
& = \int_0^{2\pi}  \int^{\beta}_{-\beta+\pi} \int_\bbR \psi(s,t,\theta')
K \big(\, \big( x+iy,e^{\frac12(y-\frac\pi2)}e^{i\theta}\big), \,
\big( s+it,e^{\frac12(t-\frac\pi2)}e^{i\theta'}\big)\, \big) \, dsdtd\theta' \\
& = \int_0^{2\pi} \int^{\beta}_{-\beta+\pi}  \int_\bbR \psi(s,t,\theta')
\Big( \sum_{j\in\bbZ} e^{(t+y-\pi)\frac j2} e^{ij(\theta-\theta')}
\frac{1}{2\pi} \int_\bbR 
\frac{e^{i(x-s)\xi  -(t+y)\xi} }{\nu(\xi,j)}\, d\xi \,
\Big) \, dsdtd\theta' \\
& = :    \int_{I_1} \int_\bbT \int_\bbR
\psi(s,t,\theta')
H(x-s,y+t,\theta-\theta') \, dsd\theta'\, dt \,.
\end{align*} 
Thus,  $T_{1,1}$ is the composition of a convolution operator on  $\bbR\times\bbT$
in the first and third variables, depending on the parameters $(y,t)$ 
 with  an integration in the $t$ variable:
$$
T_{1,1}\psi (x,y,\theta)
=  \int_{I_1}  M^{1,1}_{(y,t)} \big( \psi(\cdot, t,\cdot) \big) 
(x,\theta) 
\, dt \,,
$$
where $ M^{1,1}_{(y,t)}$ is the convolution operator on $\bbR\times\bbT$
with kernel $H(x,y+t,\theta)$.  Notice that the
Fourier multiplier associated to $M^{1,1}_{(y,t)}$, that is, the Fourier transform of  this  kernel on $\bbR\times\bbT$, is
\begin{align}
\frac{e^{(t+y-\pi)\frac j2}  e^{ -(t+y)\xi} }{\nu(\xi,j)}
&= e^{-\frac12(t+y-\pi)} 
\frac{e^{\pi\xi} e^{-\frac12(t+y-\pi)(2\xi-(j+1))} }{\nu(\xi,j)} \notag\\
& =: e^{-\frac12(t+y-\pi)} m^{1,1}_{(y,t)} (\xi,j) \,, \label{mult-m11} 
\end{align}
where $t,y\in I_1$. \ms

\noindent
{\em Remark.} We point out that, for the later part of our proof, it is convenient to
express the multiplier as function of $\xi$ and $2\xi-(j+1)$, since
the term $\nu(\xi,j)$ is in fact already function of such variables.
We will keep this approach in the expression of all the multipliers 
$m_{(y,t)}^{\ell,\ell'}$ below.  \ms
 
The operator $T_{3,3}$ is quite similar to $T_{1,1}$, and arguing in the same way we
obtain that 
\begin{align*}
T_{3,3}\psi (x,y,\theta)
& = \int_0^{2\pi}  \int^{\beta-\pi}_{-\beta} \int_\bbR \psi(s,t,\theta')
K \big(\, \big( x+iy,e^{\frac12(y+\frac\pi2)}e^{i\theta}\big), \,
\big( s+it,e^{\frac12(t+\frac\pi2)}e^{i\theta'}\big)\, \big) \,
dsdtd\theta' 
\notag \\
& = \int_0^{2\pi} \int^{\beta-\pi}_{-\beta}  \int_\bbR \psi(s,t,\theta')
\Big( \sum_{j\in\bbZ} e^{(t+y+\pi)\frac j2} e^{ij(\theta-\theta')}
\frac{1}{2\pi} \int_\bbR 
\frac{e^{i(x-s)\xi  -(t+y)\xi} }{\nu(\xi,j)}\, d\xi \,
\Big) \, dsdtd\theta' \\
& 
=  \int_{I_3} M^{3,3}_{(y,t)}\big( \psi(\cdot, t,\cdot) \big) 
(x,\theta) 
\, dt \,,
\end{align*}
where in this case $ M^{3,3}_{(y,t)}$ is the convolution operator whose associate multiplier on  $\bbR\times\bbT$ is
\begin{align}
\frac{e^{(t+y+\pi)\frac j2}  e^{ -(t+y)\xi} }{\nu(\xi,j)}
&= e^{-\frac12(t+y+\pi)} 
\frac{e^{\pi\xi} e^{-\frac12(t+y+\pi)(2\xi-(j+1))} }{\nu(\xi,j)}  \notag \\
& =: e^{-\frac12(t+y+\pi)} m^{3,3}_{(y,t)} (\xi,j) \label{mult-m33} \,,
\end{align}
with $t,y\in I_3$.
\ms

Next we turn to the case of $T_{2,2,}$.  For
$\psi \in C_c(\cX_2)$ we have that 
\begin{align*}T_{2,2}\psi (x,y,\theta)
& = \int_0^{2\pi}  \int^{\beta}_{\beta-\pi} \int_\bbR \psi(s,t,\theta')
K \big(\, \big( x+iy,e^{\frac12(\beta-\frac\pi2)}e^{i\theta}\big), \,
\big( s+it,e^{\frac12(\beta-\frac\pi2)}e^{i\theta'}\big)\, \big) \,
dsdtd\theta' 
\notag \\
& = \int_0^{2\pi} \int^{\beta}_{\beta-\pi}  \int_\bbR \psi(s,t,\theta')
\Big( \sum_{j\in\bbZ} e^{(\beta-\frac\pi2)j} e^{ij(\theta-\theta')}
\frac{1}{2\pi} \int_\bbR 
\frac{e^{i(x-s)\xi  -(t+y)\xi} }{\nu(\xi,j)}\, d\xi \,
\Big) \, dsdtd\theta' \\
& 
=  \int_{I_2}  M^{2,2}_{(y,t)}\big( \psi(\cdot, t,\cdot) \big) 
(x,\theta) 
\, dt \,,
\end{align*}
where in this case the associated multiplier is
\begin{align}
\frac{e^{(\beta-\frac\pi2)j} e^{ -(t+y)\xi} }{\nu(\xi,j)}
& = e^{-(\beta-\frac\pi2)}
\frac{e^{-(\beta-\frac\pi2)(2\xi-(j+1)) } e^{(2\beta-\pi-(t+y))\xi} }{\nu(\xi,j)}  \notag \\
& =: e^{-(\beta-\frac\pi2)} m^{2,2}_{(y,t)} (\xi,j) \label{mult-m22} \,,
\end{align}
with $t,y\in I_2$.

The case of $T_{4,4,}$ is again quite similar.  We obtain that, for
$\psi$ defined on $I_4$, 
$$
T_{4,4}\psi (x,y,\theta)
= \int_{I_4}  M^{4,4}_{(y,t)} \big( \psi(\cdot, t,\cdot) \big) 
(x,\theta) 
\, dt \,,
$$
where the Fourier multiplier of $M^{4,4}_{(y,t)}$ is
\begin{align}
\frac{e^{-(\beta-\frac\pi2)j} e^{ -(t+y)\xi} }{\nu(\xi,j)}
& = e^{(\beta-\frac\pi2)} \frac{e^{(\beta-\frac\pi2)(2\xi-(j+1)) }
  e^{-(2\beta-\pi+t+y)\xi} }{\nu(\xi,j)}
\notag\\
& =:e^{(\beta-\frac\pi2)} m^{4,4}_{(y,t)} (\xi,j)
\,,  \label{mult-m44}
\end{align}
with $t,y\in I_4$. \ms

Next we turn to the non-diagonal cases, that is, the operators
$T_{\ell,\ell'}$ with $\ell\neq\ell'$.  We proceed in the same fashion
as in the previous cases. 

We now describe $T_{2,1}$. 
For $\psi$ sufficiently regular and defined on $\cX_1$ and
$(x,y,\theta)\in\cX_2$ we have 
\begin{align*}
 T_{2,1}\psi (x,y,\theta)
& = \int_0^{2\pi}  \int_{-\beta+\pi}^\beta \int_\bbR \psi(s,t,\theta')
K \big(\, \big( x+iy,e^{\frac12(\beta-\frac\pi2)}e^{i\theta}\big), \,
\big( s+it,e^{\frac12(t-\frac\pi2)}e^{i\theta'}\big)\, \big) \,
dsdtd\theta' 
\notag \\
& =  \int_0^{2\pi}  \int_{-\beta+\pi}^\beta \int_\bbR \psi(s,t,\theta')
\Big( \sum_{j\in\bbZ} e^{(\beta-\pi+t)\frac j2} e^{ij(\theta-\theta')}
\frac{1}{2\pi} \int_\bbR 
\frac{e^{i(x-s)\xi  -(t+y)\xi} }{\nu(\xi,j)}\, d\xi \,
\Big) \, dsdtd\theta' \\
& 
=  \int_{I_1}  M^{2,1}_{(y,t)} \big( \psi(\cdot, t,\cdot) \big) 
(x,\theta) 
\, dt \,,
\end{align*}
where $ M^{2,1}_{(y,t)}$ is a convolution operator on $\bbR\times\bbT$ whose associated multiplier on $\bbR\times\bbT$ is
\begin{align}
 \frac{e^{(\beta-\pi+t)\frac j2}  e^{ -(t+y)\xi} }{\nu(\xi,j)}
& = e^{-\frac12 (\beta-\pi+t)}
\frac{e^{-\frac12 (\beta+t-\pi)(2\xi-(j+1))}e^{-(y-\beta+\pi)\xi}}{\nu(\xi,j)} \notag\\
& =: e^{-\frac12 (\beta-\pi+t)} m^{2,1}_{(y,t)}(\xi,j) \,, \label{mult-m21} 
\end{align}
where $t\in I_1$ and $y\in I_2$.

Next we focus on the operator $T_{3,1}$. We point out that the Fourier
multiplier we obtain for this operator has a better decay of the
multipliers that appeared so far. This is due to the fact that
$\operatorname{dist}(\cX_3,\cX_1)>0$. 
For $\psi$ sufficiently regular and defined on $\cX_1$ and $(x,y,\theta)\in\cX_3$ we have
\begin{align*}
 T_{3,1}\psi (x,y,\theta)
& = \int_0^{2\pi}  \int_{-\beta+\pi}^\beta \int_\bbR \psi(s,t,\theta')
K \big(\, \big( x+iy,e^{\frac12(y+\frac\pi2)}e^{i\theta}\big), \,
\big( s+it,e^{\frac12(t-\frac\pi2)}e^{i\theta'}\big)\, \big) \,
dsdtd\theta' 
\notag \\
& =  \int_0^{2\pi}  \int_{-\beta+\pi}^\beta \int_\bbR \psi(s,t,\theta')
\Big( \sum_{j\in\bbZ} e^{(y+t)\frac j2} e^{ij(\theta-\theta')}
\frac{1}{2\pi} \int_\bbR 
\frac{e^{i(x-s)\xi  -(t+y)\xi} }{\nu(\xi,j)}\, d\xi \,
\Big) \, dsdtd\theta' \\
& 
=  \int_{I_1}  M^{3,1}_{(y,t)} \big( \psi(\cdot, t,\cdot) \big) 
(x,\theta) 
\, dt \,,
\end{align*}
where the Fourier transform on $\bbR\times\bbT$ of the kernel
of the convolution operator $M^{3,1}_{(y,t)}$ is
\begin{align}
 \frac{e^{(\beta-\pi+t)\frac j2}  e^{ -(t+y)\xi} }{\nu(\xi,j)}
& = e^{-\frac12 (y+t)} \frac{e^{-\frac12 (y+t)(2\xi-(j+1))}}{\nu(\xi,j)} \notag\\
& =: e^{-\frac12 (y+t)} m^{3,1}_{(y,t)}(\xi,j) \,, \label{mult-m31} 
\end{align}
where $t\in I_1$ and $y\in I_3$.

In a similar way we obtain the operator $T_{4,1}$. 
For a sufficiently regular function $\psi$ defined
on $\cX_4$, we obtain
\begin{align*}
 T_{4,1}\psi(x,y,\theta)
= \int_{I_1} M^{4,1}_{(y,t)} \big(\psi(\cdot,t,\cdot)\big)(x,\theta)dt
\end{align*}
where the Fourier multiplier of $M^{4,1}_{(y,t)} $ is

\begin{align}
 \frac{e^{-(\beta-t)\frac j2}  e^{ -(t+y)\xi} }{\nu(\xi,j)}
& =e^{\frac12 (\beta-t)}\frac{e^{\frac12
    (\beta-t)(2\xi-(j+1))}e^{-(\beta+y)\xi}}{\nu(\xi,j)} 
\notag\\
& =: e^{\frac12 (\beta-t)} m^{4,1}_{(y,t)}(\xi,j) \,, \label{mult-m41} 
\end{align}

where $y\in I_4$ and $t\in I_1$.

All
other operators 
$T_{\ell,\ell'}$ can be obtained by slight modification of the
arguments above.

If $\psi$ is a suitable function defined on $\cX_{\ell'}$ and
$(x,y,\theta)\in\cX_\ell$, 
then a generic operator $T_{\ell,\ell'}$ turns out to be of the form
\begin{align}
T_{\ell,\ell'}\psi(x,y,\theta)=\int_{I_{\ell'}} 
M^{\ell,\ell'}_{(y,t)}\big(\psi(\cdot,t,\cdot)\big)(x,\theta)dt
\end{align}
where $M^{\ell,\ell'}_{(y,t)}$ is a Fourier multiplier operator with
multiplier of the form
$$
 \alpha^{\ell,\ell'}(y,t)\ m^{\ell,\ell'}_{(y,t)}(\xi,j)
$$
where $(y,t)\in I_\ell\times I_{\ell'}$ and the factor
$\alpha^{\ell,\ell'}$ is a positive function which
is bounded together with all its derivatives. 
Therefore, in order to prove the $L^p$ boundedness, it suffices to
consider 
the multipliers 
$m^{\ell,\ell'}_{(y,t)}$. It is now easy to see that 
\begin{align}
\begin{split}
&m^{1,2}_{(y,t)}(\xi,j)
=\frac{e^{-\frac12
    (\beta-\pi+y)(2\xi-(j+1))}e^{-(t-\beta+\pi)\xi}}{\nu(\xi,j)}\,,
\qquad 
 m^{1,4}_{(y,t)}(\xi,j)=\frac{e^{\frac12
    (\beta-y)(2\xi-(j+1))}e^{-(\beta+t)\xi}}{\nu(\xi,j)}\,,\\ 
&m^{3,2}_{(y,t)}(\xi,j)
=\frac{e^{-\frac12
    (\beta+y)(2\xi-(j+1))}e^{(\beta-t)\xi}}{\nu(\xi,j)}\,, 
\qquad
m^{2,3}_{(y,t)}(\xi,j)=\frac{e^{\frac12
    (\beta+t)(2\xi-(j+1))}e^{(\beta-y)\xi}}{\nu(\xi,j)}\,,\\ 
&m^{3,4}_{(y,t)}(\xi,j)
=\frac{e^{\frac12
    (\beta-\pi-y)(2\xi-(j+1))}e^{-(\beta-\pi+t)\xi}}{\nu(\xi,j)} \,,
 \quad 
m^{4,3}_{(y,t)}(\xi,j)=\frac{e^{\frac12
     (\beta-\pi-t)(2\xi-(j+1))}e^{-(\beta-\pi+y)\xi}}{\nu(\xi,j)} \,, 
 \label{multipliers}
\\
 & m^{1,3}_{(y,t)}(\xi,j)=\frac{e^{-\frac12(t+y)(2\xi-(j+1))}}{\nu(\xi,j)}
\,, \qquad
m^{4,2}_{(y,t)}(\xi,j)=\frac{e^{(y+t)\xi}}{\nu(\xi,j)}\,, \qquad 
 m^{2,4}_{(y,t)}(\xi,j)=\frac{e^{-(y+t)\xi}}{\nu(\xi,j)} \,.
 \end{split}
\end{align}

{\bf (3)}
We proceed by applying Minkowski's integral inequality.
Let $\ell,\ell'\in\{1,\dots,4\}$. Then
$$
|T_{\ell,\ell'}\psi (x,y,\theta)| \le
 \int_{I_{\ell'}} \big|M^{\ell,\ell'}_{(y,t)}\big( \psi(\cdot, t,\cdot) \big) 
(x,\theta) \big|
\, dt \,,
$$
so that
\begin{align*}
\int_{\cX_\ell}   |M^{\ell,\ell'}_{(y,t)}\psi (x,y,\theta)|^p \, dxdyd\theta 
& \le   \int_{I_\ell}  \int_{\bbR\times\bbT} \bigg(
 \int_{I_{\ell'}}\big|M^{\ell,\ell'}_{(y,t)} \big( \psi(\cdot, t,\cdot) \big) 
(x,\theta) \big|
\, dt \bigg)^p \, dxd\theta \,dy  \\
& =  \int_{I_\ell}   \bigg\{ \int_{\bbR\times\bbT} \bigg(
 \int_{I_{\ell'}}\big|M^{\ell,\ell'}_{(y,t)}\big( \psi(\cdot, t,\cdot) \big) 
(x,\theta) \big|
\, dt \bigg)^p \, dxd\theta \bigg\}^{p/p}\,  dy \\
& \le \int_{I_\ell} \bigg\{ 
\int_{I_{\ell'}} \bigg(  \int_{\bbR\times\bbT} \big|M^{\ell,\ell'}_{(y,t)}\big( \psi(\cdot, t,\cdot) \big) 
(x,\theta) \big|^p
\, dxd\theta \bigg)^{1/p} \,dt
\bigg\}^p \, dy 
\,.
\end{align*}

Then, if we denote by $N_p^{\ell,\ell'}(y,t)$ the operator norm of $M^{\ell,\ell'}_{(y,t)}$ on
$L^p(\bbR\times\bbT)$ we have
\begin{align}
\| T_{\ell,\ell'}\psi\|_{L^p(\cX_{\ell})}^p
& = \int_{\cX_\ell}  |T_{\ell,\ell'}\psi (x,y,\theta)|^p \, dxdyd\theta \notag\\
& \le \int_{I_\ell} \bigg\{ 
\int_{I_{\ell'}} N_p^{\ell,\ell'}(y,t) \| \psi(\cdot, t,\cdot) \|_{L^p(\bbR\times\bbT)}   \,dt
\bigg\}^p \, dy  \,.
\label{Main Estimate}
\end{align}

Thus, we have reduced ourselves to showing that:
\begin{itemize} 
\item[{\bf (A)}]
the operator $M^{\ell,\ell'}_{(y,t)}$ is bounded on $L^p(\bbR\times\bbT)$ for
a.a. $t\in I_{\ell'}$ and $y\in I_\ell$
with norm $N_p^{\ell,\ell'}(y,t)$;
\item[{\bf (B)}] the integral operator $\cN_{\ell,\ell'}$ with
  integral kernel
  $N_p^{\ell,\ell'}(y,t)$,  is bounded  
 $\cN_{\ell,\ell'}: L^p(I_{\ell'}) \to L^p(I_\ell)$,   for $1<p<\infty$. \ms
 \end{itemize}

>From these two facts, the boundedness of $T_{\ell,\ell'}$, hence the one of
$P$, follows at once.  

We point that the integral operator $\cN_{\ell,\ell'}$  with integral kernel
  $N_p^{\ell,\ell'}(y,t)$ turns to be an operator of Hilbert type
  (see \cite{Phong-Stein}).  This fact is not surprising as several
  operators in several complex variables appear as composition of
  singular integrals (as the operator $M^{\ell,\ell'}_{(y,t)}$) and operators of
  Hilbert type.
\ms

In order to prove {\bf (A)}, recall that $M^{\ell,\ell'}_{(y,t)}$ is a
Fourier multiplier operator on $\bbR\times\bbZ$, 
whose multiplier is given in the identities 
(\ref{mult-m11}\,-\,\ref{mult-m41}) and \eqref{multipliers}.
By the boundedness of the factors appearing on the right hand sides of
such equations, it
 suffices to study the boundedness of the
multipliers $m_{(y,t)}^{\ell,\ell'}(\xi,j)$, for $\ell,\ell'\in
\{1,\dots,4\}$. 

 It suffices
to show that $m_{(y,t)}^{\ell,\ell'}(\xi,\eta)$ is a bounded Fourier multiplier on
$L^p(\bbR^2)$, see \cite[Thm. 3.6.7]{MR2445437}.  
By the  affine change of variables 
\begin{equation}\label{aff-chg-coor}
\begin{cases}
\xi'= \pi\xi  \cr
\eta'=(\beta-\pt)(2\xi-(\eta+1))\,,
\end{cases}
\end{equation} 
 it suffices to show
that the functions 
$\tilde m_{(y,t)}^{\ell,\ell'} (\xi',\eta')$ 
give rise to bounded Fourier multipliers on $L^p(\bbR^2)$.

\ms

Denote by $T_m$ the Fourier multiplier operator with multiplier
$m$.  Then, using \cite[Thm. 3.6.7]{MR2445437} again, we have the bounds
$$
\| M^{\ell,\ell'}_{(y,t)} \|_{L^p(\bbR\times\bbZ)\to L^p(\bbR\times\bbZ)}
\le \| T_{m_{(y,t)}^{\ell,\ell'}} \|_{L^p(\bbR^2)\to L^p(\bbR^2)}
\le C \| T_{\tilde m_{(y,t)}^{\ell,\ell'}} \|_{L^p(\bbR^2)\to
  L^p(\bbR^2)} \,,
$$
with $C$ independent of $y$ and $t$.   \ms

\noindent
{\em Proof of } {\bf (A)}.
We shall show that $\tilde m^{\ell,\ell'}_{(y,t)}$ is a Marcinkiewicz multiplier in
$\bbR^2$ with norm $N_p^{\ell,\ell'}(y,t)$, for $\ell,\ell'\in\{1,\dots,4\}$.

We recall that  $m$ is said to be a  
Marcinkiewicz multiplier in $\bbR^2$ if there exists $C>0$ such that 
\begin{equation}\label{Mar-est}
\sup_{\xi,\eta\in\bbR\setminus\{0\}}  \big|\p^{k_1}_\xi \p^{k_2}_\eta  
m (\xi,\eta)\big|  \le C \frac{1}{ |\xi|^{k_1} |\eta|^{k_2} } 
\end{equation}

for $k_1,k_2\ge0$, $k_1+k_2\le2$, see
e.g. \cite[Corollary 5.2.5]{MR2445437}.  Thus, we wish to show that 
$\tilde m_{(y,t)}^{\ell,\ell'}$ are Marcinkiewicz multipliers in
$\bbR^2$ and estimates their operator norms, for
$\ell,\ell'\in\{1,\dots,4\}$. 
\ms

In view of the change of
coordinates \eqref{aff-chg-coor} (and dropping the primes to  ease 
notation), we 
 set
\begin{equation}\label{Denominator}
D(\xi,\eta)
=\cosh\eta\frac{\sinh\xi}{\xi}+\frac{\sinh\eta}{\eta}\cosh\xi \,,
\end{equation}
and
\begin{equation*}
Q(\xi,\eta)
=\frac{D(\xi,\eta)}{\tilde\nu(\xi,\eta)}\,,
\end{equation*}
where $\tilde\nu(\xi,\eta)$ denotes the function $\nu(\xi,\eta)$ 
after the change of variables \eqref{aff-chg-coor}.

In the diagonal cases, the multipliers 
$\tilde m_{(y,t)}^{1,1},\dots,\tilde m_{(y,t)}^{4,4}$, are given by
(\ref{mult-m11}\,-\,\ref{mult-m44}), namely,
\begin{equation*}
\tilde m_{(y,t)}^{\ell,\ell}(\xi,\eta)
= \frac{e^{\gamma\xi} e^{\alpha\eta} }{D(\xi,\eta)}\cdot Q(\xi,\eta) \,,
\end{equation*}
 where
\begin{align}\label{Tilde Symmetric}
\begin{split}
& \text{if } \ell=1\,, \qquad \alpha =\frac{\pi-(t+y)}{2\beta-\pi}\,, \gamma=1\,, \ t,y\in I_1 \,;\\
& \text{if } \ell=2\,, \qquad \alpha=-1\,, \gamma =\frac{2\beta-\pi-(t+y)}{\pi}
\,, 
\ t,y\in I_2
\,;\\
& \text{if } \ell=3\,, \qquad \alpha =-\frac{t+y+\pi}{2\beta-\pi}\,, \gamma=1\,, \ t,y\in I_3 \,;\\
& \text{if } \ell=4\,, \qquad \alpha=1\,, \gamma =-\frac{2\beta-\pi+ t+y}{\pi}
\,, 
\ t,y\in I_4
\,.
\end{split}
\end{align}

Next we turn to the off-diagonal cases. As in the diagonal case, from
(\ref{mult-m11}\,-\,\ref{mult-m41}),  \eqref{multipliers} and the
change of variables 
\eqref{aff-chg-coor}, we obtain again
\begin{align*}
 \tilde m^{\ell,\ell'}_{(y,t)}(\xi,\eta)
=\frac{e^{\gamma\xi}e^{\alpha\eta}}{D(\xi,\eta)}\cdot Q(\xi,\eta),
\end{align*}
where now,
\begin{align*}
\begin{split}
 &\text{if } (\ell,\ell')=(2,1)\,, \qquad \alpha =-\frac{\beta-\pi+t}{2\beta-\pi}\,, \gamma=-\frac{y-\beta+\pi}{\pi}\,, \ t\in I_1 \,, y\in I_2 \,;\\
 &\text{if } (\ell,\ell')=(4,1)\,, \qquad \alpha =-\frac{\beta-t}{2\beta-\pi}\,, \gamma=-\frac{\beta+y}{\pi}\,, \ t\in I_1 \,, y\in I_4 \,;\\
 &\text{if } (\ell,\ell')=(1,2)\,, \qquad \alpha =-\frac{\beta-\pi+y}{2\beta-\pi}\,, \gamma=-\frac{t-\beta+\pi}{\pi}\,, \ t\in I_2 \,, y\in I_1 \,;\\
 &\text{if } (\ell,\ell')=(3,2)\,, \qquad \alpha =-\frac{\beta+y}{2\beta-\pi}\,, \gamma=\frac{\beta-t}{\pi}\,, \ t\in I_2 \,, y\in I_3 \,;\\
 &\text{if } (\ell,\ell')=(2,3)\,, \qquad \alpha =\frac{\beta+t}{2\beta-\pi}\,, \gamma=-\frac{\beta-y}{\pi}\,, \ t\in I_3 \,, y\in I_2 \,;\\
 &\text{if } (\ell,\ell')=(4,3)\,, \qquad \alpha =\frac{\beta-\p-t}{2\beta-\pi}\,, \gamma=\frac{\pi-\beta-y}{\pi}\,, \ t\in I_3 \,, y\in I_4 \,;\\
 &\text{if } (\ell,\ell')=(1,4)\,, \qquad \alpha =\frac{\beta-y}{2\beta-\pi}\,, \gamma=-\frac{\beta+t}{\pi}\,, \ t\in I_4 \,, y\in I_1 \,;\\
 &\text{if } (\ell,\ell')=(3,4)\,, \qquad \alpha =\frac{\beta-\pi-y}{2\beta-\pi}\,, \gamma=\frac{\pi-\beta-t}{\pi}\,, \ t\in I_4 \,, y\in I_3 \,.\\
 \end{split}
\end{align*}

We point out that the cases when $(\ell,\ell')\in
\{(3,1),(4,2),(1,3),(2,4)\}$ are easier and the
multipliers have a better decay. This is due to the fact that, as we
already mentioned, 
$\operatorname{dist}(\cX_1,\cX_3)$,
$\operatorname{dist}(\cX_2,\cX_4)>0$. 
In these special cases we have
\begin{align*}
 \tilde m^{\ell,\ell'}_{(y,t)}(\xi,\eta)
=\frac{e^{\gamma\xi}}{D(\xi,\eta)}\cdot Q(\xi,\eta),
\end{align*}
where
\begin{align*}
 \begin{split}
 &\text{if } (\ell,\ell')=(3,1)\,, \qquad \gamma=-\frac{t+y}{2\beta-\pi}\,, \ t\in I_1 \,, y\in I_3 \,.\\  
 &\text{if } (\ell,\ell')=(4,2)\,, \qquad \gamma=\frac{t+y}{\pi}\,, \ t\in I_2 \,, y\in I_4 \,.\\
 &\text{if } (\ell,\ell')=(1,3)\,, \qquad \gamma=-\frac{t+y}{2\beta-\pi}\,, \ t\in I_3 \,, y\in I_1 \,.\\ 
 &\text{if } (\ell,\ell')=(2,4)\,, \qquad \gamma=-\frac{t+y}{\pi}\,, \ t\in I_4 \,, y\in I_3 \,. 
 \end{split}
\end{align*}

We shall make use of the following
lemmas and defer their proofs to Section \ref{Lemmas-Sec}. 
\begin{lem}\label{Mar-lem-Q}
 The function 
 \begin{align*}
   Q(\xi,\eta)=\frac{\cosh\eta\frac{\sinh\xi}{\xi}+\frac{\sinh\eta}{\eta}\cosh\xi}{2\pi\frac{\sinh\xi}{\xi}
\cosh\big[\eta+\frac12(\beta-\frac\pi2)\big]+(\beta-\frac\pi2)\cosh\xi
\Big(\int_{-1}^1 e^{-\eta s}\sqrt{1+e^{-(\beta-\frac\pi2)s}ds} \Big)}
  \end{align*}
is a Marcinkiewicz multiplier in $\bbR^2$.
\end{lem}

\begin{lem}\label{Mar-lem-diag}
The function
\begin{equation}\label{malpha} 
 m_\alpha
 (\xi,\eta)=\frac{e^{\xi}e^{\alpha\eta}}{D(\xi,\eta)}
\end{equation}
is a Marcinkiewicz multiplier in $\mathbb{R}^2$ for every
$0<\alpha<1$. 
Moreover the operator norm of the associated multipliers operator is
bounded 
by $\frac{C}{1-\alpha}$.
\end{lem}

\begin{lem}\label{Mar-lem-off-diag}
The function 
 \begin{equation*}
  m_{\alpha,\gamma}(\xi,\eta)
=\frac{e^{\gamma\xi}e^{\alpha\eta}}{D(\xi,\eta)}
 \end{equation*}
is a Marcinkiewicz multiplier on $\mathbb{R}^2$ for every
$0<\alpha,\gamma<1$. Moreover the norm of the associated multiplier
operator is bounded by $\frac{C}{(1-\alpha)+(1-\gamma)}$. 
\end{lem}

 \ms

\noindent
{\em Proof of } {\bf (B)}. We now prove {\bf (B)}
in the case of the operator $\cN_{1,1}$
associated to $T_{1,1}$. The remaining cases can be proved similarly
with minor modifications.  
We show that  $\cN_{1,1}$ is a Hilbert-type integral operator
and can be studied with standard technique by means of Schur's test
(see, e.g., \cite[Appendix I]{MR2445437} ). 

 We suppose that the parameter
$\alpha$ appearing in $\tilde m^{1,1}_{y,t}$, see \eqref{Tilde
  Symmetric}, is positive. The case of $\alpha$ negative can be
similarly treated. From \eqref{Main Estimate}, \eqref{Tilde Symmetric}
and Lemma \ref{Mar-lem-diag} we deduce 
\begin{align*}
 \| T_{1,1}\psi\|^p_{L^p(\cX_\ell}&\leq C \int_{I_1}\bigg\{\int_{I_1}
 \frac{1}{1-\frac{\pi-(t+y)}{2\beta-\pi}}\|\psi(\cdot,t,\cdot)\|_{L^{p}(\bbR\times\bbT)}\
 dt\bigg\}^p\ dy\\ 
 &\leq
 C\int_{I_1}\bigg\{\int_{I_1}\frac{1}{(\beta-\pi+t)+(\beta-\pi+y)}\|
\psi(\cdot,t,\cdot)\|_{L^p(\bbR\times\bbT)}\  dt\bigg\}^p \
 dy. 
\end{align*}
Therefore, the integral operator $\cN_{1,1}$ associated to $T_{1,1}$ has
positive integral kernel  
$$
N_p^{1,1}(y,t)=\frac{1}{(\beta-\pi+t)+(\beta-\pi+y)}\,.
$$
It suffices to show that there
exists a constant $C>0$ such that 
\begin{align*}
 \int_{I_1}\frac{1}{(\beta-\pi+t)+(\beta-\pi+y)} 
\varphi^p(t)\ dt\leq C\varphi^p(y)
\end{align*}
and 
\begin{align*}
 \int_{I_1}\frac{1}{(\beta-\pi+t)+(\beta-\pi+y)}
\varphi^{q}(t)\ dt\leq C \varphi^q(y)
\end{align*}
for every $y\in I_1$, where $q$ is the conjugate exponent of $p$. It
is enough to choose $\varphi(t)=(\beta-\pi+t)^{-\frac{1}{qp}}$ and the
desired estimates are easily proved; see \cite[(2.2)]{Phong-Stein} 
or
\cite[Appendix I.3]{MR2445437} for details.

 The
estimates for the operators $T_{\ell,\ell}$, $\ell=2,3,4$ follow in
the same fashion.

When $\ell\neq\ell$,  it suffices to
perform an affine change of variables so that $I_{\ell}$ coincides
with $I_{\ell'}$ and to use again Schur's test.  We leave the
elementary details to the reader.
\qed

\ms

\section{Proofs of Theorem \ref{Sze-pro-main-thm} and
 Lemmas \ref{Mar-lem-Q}-\ref{Mar-lem-off-diag}} 
\label{Lemmas-Sec}
  \ms

\subsection{Proof of  Theorem \ref{Sze-pro-main-thm}}
In order to complete the proof it suffices to show that 
$$
P: W^{k,p}(bD'_\beta)\to
 W^{k,p}(bD'_\beta)
$$
 is bounded when
$1<p<\infty$ and $k$ is a positive integer.  Clearly, again it also suffices
to show that the operators $T_{\ell,\ell'}$ in Theorem 
\ref{Lp-bddness} are bounded on $W^{k,p}(bD'_\beta)$, when
$1<p<\infty$ and $k$ is a positive integer.

Since $T_{\ell,\ell'} = \int_{I_{\ell'}} M^{\ell,\ell'}_{(y,t)} \,dt$,
where $M^{\ell,\ell'}_{(y,t)}$ is a Fourier multiplier operator on
$\bbR\times\bbT$, derivatives in $x$ and $\theta$ commute with
$T_{\ell,\ell'}$ and therefore we have nothing to prove.  

Finally, consider derivatives in $y\in I_\ell$.  
We claim that $\p_y (T_{\ell,\ell'} f)= T_{\ell,\ell'}(\p_x f)$, up to
a multiplicative constant.  This commutation rule will prove the theorem.

For,
when we examine the
dependence of the multiplier $m_{(t,y)}^{\ell,\ell'}$ in $y$ we see
that they either have no effect on the multiplier itself (if the
derivative falls on a factor of the form $e^{cy}$, $c$ a constant), or
produce a term that has the form $\xi \cdot m_{(t,y)}^{\ell,\ell'}$
(when the derivative falls on a factor of the form $e^{y\xi}$).
Multiplication by $\xi$ on the Fourier transform side can be viewed as
a derivative in $x$ of the function the operator $T_{\ell,\ell'}$ is
acting on, and we are done.  \qed
\ms

\subsection{Proof of Lemma \ref{Mar-lem-Q}}
It is easy to see that
$$
\frac{D(\xi,\eta)}{2(\beta-\frac\pi2)\cosh\eta\frac{\sinh\eta}{\eta}+2\pi\frac{\sinh\xi}{\xi}\cosh\eta}
$$
is a Marcinkiewicz
multiplier, therefore, in order to prove that $Q(\xi,\eta)$ is
a multiplier, it suffices to focus on
\begin{align*}
 \frac{2(\beta-\frac\pi2)\cosh\eta\frac{\sinh\eta}{\eta}
+2\pi\frac{\sinh\xi}{\xi}\cosh\eta}{\tilde\nu(\xi,\eta)}
  & =\frac{2(\beta-\frac\pi2)
   \frac{\tanh\eta}{\eta}+2\pi\frac{\tanh\xi}{\xi}}{2\pi\frac{\tanh\xi}{\xi}
   +(\beta-\frac\pi2)F(\eta)}\cdot 
\frac{\cosh\eta}{\cosh(\eta+\frac12(\beta-\frac\pi2))}\\
 & =:\tilde Q(\xi,\eta)\cdot P(\eta),
\end{align*}
 where
\begin{align*}
F(\eta)
&=\frac{\int_{|s|<1}e^{-\eta s}\sqrt{1+4^{-(\beta-\frac\pi2)s}}\, ds}{ \cosh(\eta+\frac12(\beta-\frac\pi2))}\\
&=\frac{\int_{|s|<1}e^{-\eta s}\sqrt{1+4^{-(\beta-\frac\pi2)s}}\,
  ds}{\cosh\eta} \cdot
\frac{\cosh\eta}{\cosh(\eta+\frac12(\beta-\frac\pi2))}\\
&=: \tilde F(\eta)\cdot P(\eta).
\end{align*}
The factor $P(\eta)$ in the above formulas is clearly bounded and it
is easily verified that all its derivatives decay
exponentially. 

In order to show that $\tilde Q(\xi,\eta)$ is a
Marcinkiewicz multiplier, we write
\begin{align*}
 \tilde Q(\xi,\eta)
&=1+(\beta-\frac\pi2) 
\bigg[ 
\frac{2\frac{\tanh\eta}{\eta}}{2\pi\frac{\tanh\xi}{\xi} +(\beta-\frac\pi2)F(\eta)} 
-\frac{F(\eta)}{2\pi\frac{\tanh\xi}{\xi}+(\beta-\frac\pi2)F(\eta)} 
\bigg]\\
 &=: 1+(\beta-\frac\pi2)\big[\tilde Q_1(\xi,\eta)+\tilde Q_2(\xi,\eta)\big].
\end{align*}
Therefore, we reduced our problem to prove that the functions $\tilde
Q_1$ and $\tilde Q_2$ are Marcinkiewicz multipliers.  

We now explicitly study the decay of the derivatives of the function $F(\eta)$. Once the behavior of $F(\eta)$ is known, it is straightforward to prove that $\tilde Q$, hence $Q$, is a multiplier.

It is immediately seen that there exist two positive constants $c_1,c_2$ such that 
$$
c_1\frac{\tanh\eta}{\eta}\leq F(\eta)\leq c_2 \frac{\tanh\eta}{\eta}
\,.
$$
Next, we focus on the derivative of $F(\eta)$ and we  may
just consider the factor
$\tilde F(\eta)$ above.

We 
assume $\eta>0$, the case $\eta<0$ can be treated in the same fashion.
We write
\begin{align*}
\tilde F(\eta)&= \frac{\int_0^1 e^{-\eta
    s}\sqrt{1+4e^{-(\beta-\frac\pi2)s}}\ ds+\int_0^1 e^{\eta
    s}\sqrt{1+4e^{(\beta-\frac\pi2)s}}\ ds}{\cosh\eta}\\ 
 &=\frac{2}{1+e^{-2\eta}}\bigg(\int_{0}^1 e^{-\eta(s+1)}\sqrt
 {1+4^{-(\beta-\frac\pi2)s}}\  ds+\int_{0}^1
 e^{\eta(s-1)}\sqrt{1+4e^{(\beta-\frac\pi2)s}}\ ds\bigg)\\ 
 &=: I(\eta)+I\!I(\eta).
\end{align*}
Now, $\frac{dI}{d\eta}(\eta)= A(\eta)-B(\eta)$, where
\begin{align*}
&A(\eta)=\frac{4e^{-2\eta}}{(1+e^{-2\eta})^2} 
\int_0^1 e^{-\eta(s+1)}\sqrt{1+4e^{(-\beta-\frac\pi2)s}}\ ds;\\
&B(\eta)=\frac{2}{1+e^{-2\eta}} 
\int_0^1 (s+1)e^{-\eta(s+1)}\sqrt{1+ 4e^{-(\beta-\frac\pi2)s}}\ ds.
\end{align*}
The term $A(\eta)$ decays exponentially, while 
\begin{align*}
 |B(\eta)|
&\leq \frac{1}{1+e^{-2\eta}}\int_0^2 \tau
e^{-\eta\tau}\sqrt{1+4e^{-(\beta-\frac\pi2)(\tau-1)}}\ d\tau\\ 
 &\leq \frac{C}{\eta^2} \int_0^{2\eta} t e^{-t}dt\\
 &\leq \frac{C}{\eta^2}.
\end{align*}
Therefore, we conclude that
$|\frac{dI}{d\eta}(\eta)|\leq\frac{C}{\eta^2}$. Similarly, we obtain the same estimate for
$ \frac{dI\!I}{d\eta}(\eta)$, so  we conclude that
$|\tilde F' (\eta)|\leq \frac{C}{\eta^2}$. 

Arguing as above, the estimate $|\tilde F''(\eta)|\leq \frac{C}{|\eta|^3}$ can be obtained, but we leave the details to the reader. 

Finally, using the estimates we have on $\tilde F$ and its derivatives, it is now straightforward to prove that $\tilde Q$, thus $Q$, is a Marcinkiewicz multiplier and conclude the proof.
\qed

\ms

\subsection{Proof of Lemma \ref{Mar-lem-diag} }
We need to show that $m_\alpha$ given by \eqref{malpha} is a
Marcinkiewicz multiplier in $\bbR^2$.  Notice that since $\alpha>0$,
in the variable
$\eta$
we only need to
estimate the decay for $\eta>0$.   For, let $\varrho$ be a
non-negative $C^\infty$ function on $\bbR$, $\varrho= 1$ for
$\eta\ge1$ and $\varrho=0$ for $\eta\le0$.  Then $(1-\varrho)m_\alpha$
is a Schwartz function in $\bbR^2$ uniformly in $\alpha$, and we only
need to consider $\varrho(\eta)\cdot m_\alpha(\xi,\eta)$.  However, 
$\varrho'$ is compactly supported, so for simplicity of notation we
just consider $m_\alpha$ and $\eta\ge1$.

We now notice
that for every  non-negative  integer $k$ there exists
$C=C_k>0$ such that
\begin{align}\label{malphasup}
|m_\alpha(\xi,\eta)|
&\leq \frac{\eta e^{\alpha\eta}}{\sinh\eta}
= 
\left|\frac{\eta^{k+1}e^{\alpha\eta}}{\sinh\eta}\right|
\frac{1}{|\eta|^k} 
\leq \frac{C_k}{(1-\alpha)^{k+1}}\frac{1}{|\eta|^k} \,.
 \end{align}
Hence, in order to prove our lemma we need to estimate all the derivatives up
to order two of 
$m_{\alpha}(\xi,\eta)$. 
Recall that $D(\xi,\eta)$ is defined in 
\eqref{Denominator}.
We have 
\begin{align}
 \p_\xi m_\alpha(\xi,\eta)
&=e^\xi e^{\alpha\eta}\frac{D(\xi,\eta)-\p_\xi
  D(\xi,\eta)}{D^2(\xi,\eta)} \notag \\
 &=m_\alpha(\xi,\eta) \left[1-\frac{\p_\xi
     D(\xi,\eta)}{D(\xi,\eta)}\right] \notag \\
 &=
 \frac{m_\alpha(\xi,\eta)}{D(\xi,\eta)}
\left[\frac{\sinh\eta}{\eta}e^{-\xi}-\cosh\eta\frac{e^{-\xi}}{\xi}
+\cosh\eta\frac{\sinh\xi}{\xi^2}\right]  \notag \\
 &=\frac{e^{\alpha\eta}}{\big[\cosh\eta\frac{\sinh\xi}{\xi}
+\frac{\sinh\eta}{\eta}\cosh\xi\big]^2}\left[\frac{\sinh\eta}{\eta}
-\frac{\cosh\eta}{\xi}+\cosh\eta\frac{\sinh\xi}{\xi^2}e^{\xi}\right] \notag \\
 &=\psi_\alpha (\xi,\eta)
\left[\frac{\sinh\eta}{\eta}
-\frac{\cosh\eta}{\xi}+\cosh\eta\frac{\sinh\xi}{\xi^2}e^\xi
 \right] \,,\label{Derivative-xi}
 \end{align}
 where 
$$
\psi_\alpha(\xi,\eta)=\frac{e^{\alpha\eta}}{D^2(\xi,\eta)} \,. 
$$

Hence,  recalling \eqref{Denominator} and using the Cauchy-Schwarz inequality we have, 
\begin{align}
|\p_\xi m_\alpha(\xi,\eta)|
& \le \frac{e^{\alpha\eta}}{D(\xi,\eta)} \cdot
\frac{\big| \frac{\sinh\eta}{\eta}
-\frac{\cosh\eta}{\xi}+\cosh\eta\frac{\sinh\xi}{\xi^2}e^\xi
\big|}{D(\xi,\eta)}  \notag \\
&  =  \frac{\eta e^{\alpha\eta}}{\sinh\eta} \cdot\frac{1}{|\xi|}\cdot
\frac{\frac{\sinh\eta}{\eta} \big| \xi\frac{\sinh\eta}{\eta}
-\cosh\eta +\cosh\eta\frac{\sinh\xi}{\xi}e^\xi
\big|}{D(\xi,\eta)^2}
 \notag\\
& 
\le  \frac{C}{1-\alpha}\cdot\frac{1}{|\xi|} \label{der-xi}\,.
\end{align}
 
For later reference, we  observe that from \eqref{malphasup} (with
$k=1$) we also have
\begin{equation}\label{later-ref}
|\p_\xi m_\alpha(\xi,\eta)|
\le  \frac{C}{(1-\alpha)^2}\frac{1}{|\xi\eta|} \,.
\end{equation}

Next, 
\begin{align*}
& \p^2_\xi m_\alpha (\xi,\eta)\\
&= \p_\xi \psi_\alpha(\xi,\eta)\left[\frac{\sinh\eta}{\eta}
-\frac{\cosh\eta}{\xi}+\cosh\eta\frac{\sinh\xi}{\xi^2}  e^\xi  \right] \notag \\
&\qquad\qquad + \psi_\alpha(\xi,\eta)
\left[
\frac{\cosh\eta}{\xi^2}+\cosh\eta\frac{\cosh\xi}{\xi^2} e^\xi 
-2\cosh\eta\frac{\sinh\xi}{\xi^3} e^\xi 
+  \cosh\eta\frac{\sinh\xi}{\xi^2}e^\xi
\right] \notag \\
& = \psi_\alpha(\xi,\eta) \bigg\{  
\frac{-2 \big[\cosh\eta\big(\frac{\cosh\xi}{\xi}-\frac{\sinh\xi}{\xi^2}\big)
+ \frac{\sinh\eta}{\eta}\sinh\xi
\big] }{D(\xi,\eta)} \left[\frac{\sinh\eta}{\eta}
-\frac{\cosh\eta}{\xi}+\cosh\eta\frac{\sinh\xi}{\xi^2}  e^\xi 
\right]   \notag \\
&\qquad\qquad \qquad\qquad 
+ \left[
\frac{\cosh\eta}{\xi^2}+\cosh\eta\frac{\cosh\xi}{\xi^2} e^\xi 
-2\cosh\eta\frac{\sinh\xi}{\xi^3} e^\xi 
+  \cosh\eta\frac{\sinh\xi}{\xi^2}e^\xi
\right]
\bigg\} \notag\\
& = \frac{\eta e^{\alpha\eta}}{\sinh\eta} \cdot\frac{1}{\xi^2}\cdot
\frac{\frac{\sinh\eta}{\eta} }{D^3(\xi,\eta)}
\bigg\{
-2\left[
\cosh\eta\big(\cosh\xi-\frac{\sinh\xi}{\xi} \big)
+ \frac{\sinh\eta}{\eta}\xi \sinh\xi
\right]\\
&\qquad\qquad \qquad\qquad  \quad\quad\times\left[\xi\frac{\sinh\eta}{\eta}
-\cosh\eta+\cosh\eta\frac{\sinh\xi}{\xi}e^\xi\right] \notag \\
&\qquad\qquad \qquad\qquad +
\left[ \cosh\xi\frac{\sinh\eta}{\eta}+ \cosh\eta\frac{\sinh\xi}{\xi} \right]\\
&\qquad\qquad \qquad\qquad  \quad\quad\times
\left[
\cosh\eta +\cosh\eta \cosh\xi  e^\xi 
-2\cosh\eta\frac{\sinh\xi}{\xi} e^\xi 
+  \cosh\eta \sinh\xi e^\xi
\right] 
\bigg\} \\
& =  \frac{\eta e^{\alpha\eta}}{\sinh\eta} \cdot\frac{1}{\xi^2}\cdot
\frac{2\frac{\sinh\eta}{\eta} }{D^3(\xi,\eta)} 
\bigg\{ \cosh\eta\frac{\sinh\eta}{\eta} \xi e^{-\xi} +\cosh^2\eta \cosh\xi
+ \cosh\eta\frac{\sinh\eta}{\eta}\big( e^\xi+\sinh\xi\big) \\
&\qquad\qquad \qquad\qquad  -
\cosh\xi\cosh\eta\frac{\sinh\xi}{\xi}\frac{\sinh\eta}{\eta} e^\xi
- \cosh^2\eta\frac{\sinh\xi}{\xi}  \bigg\} 
\,.
 \end{align*}
Then, by \eqref{malphasup} (with $k=0$) we have
$$
|\p^2_\xi m_\alpha(\xi,\eta)|\leq \frac{C}{1-\alpha}\frac{1}{|\xi|^2} \,.
$$
\ms

Next,
\begin{align}
 \p_\eta
m_\alpha(\xi,\eta)
&=\frac{m_\alpha(\xi,\eta)}{D(\xi,\eta)}
\left[
\frac{\sinh\xi}{\xi}(\alpha\cosh\eta-\sinh\eta)
+\frac{\cosh\xi}{\eta}(\alpha\sinh\eta-\cosh\eta) 
+\cosh\xi\frac{\sinh\eta}{\eta^2}
\right] 
\notag \\
&=\frac{m_\alpha(\xi,\eta)}{D(\xi,\eta)}\bigg[(\alpha-1)D(\xi,\eta)
+e^{-\eta}\Big(\frac{\sinh\xi}{\xi}-\frac{\cosh\xi}{\eta}\Big)
+\cosh\xi \frac{\sinh\eta}{\eta^2}\bigg] \notag \\
&=(\alpha-1)m_\alpha(\xi,\eta)
+\frac{m_\alpha(\xi,\eta)}{D(\xi,\eta)}
\bigg[e^{-\eta}\Big(\frac{\sinh\xi}{\xi}
-\frac{\cosh\xi}{\eta}\Big)+\cosh\xi\frac{\sinh\eta}{\eta^2}\bigg]\notag \\ 
&=: 
\left[ (\alpha-1)+ \frac{\phi(\xi,\eta)}{D(\xi,\eta)} \right]  m_\alpha(\xi,\eta) \,.
\label{Derivative-eta}
\end{align}
Using \eqref{malphasup} (with $k=0,1$) and the fact that $\eta\ge1$ so that
$|\phi/D|\le C|\eta|^{-1}$, we obtain
$$
|\p_\eta m_\alpha(\xi,\eta)|\leq\frac{C}{1-\alpha}\frac{1}{|\eta|} \,.
$$

Furthermore,
\begin{align}
\p^2_\eta m_\alpha(\xi,\eta)
&=\!\left[ (\alpha-1)+ \frac{\phi(\xi,\eta)}{D(\xi,\eta)} \right] \p_\eta m_\alpha(\xi,\eta)
+\frac{m_\alpha(\xi,\eta)}{D^2(\xi,\eta)}\big[ \p_\eta  
\phi(\xi,\eta)D(\xi,\eta)-\phi(\xi,\eta)\p_\eta D(\xi,\eta)\big] \notag \\
&=: I+I\!I \,.\label{2-Derivative-eta}
\end{align}
From \eqref{Derivative-eta} it follows that
\begin{align*} \begin{split}
I&= \left[ (\alpha-1)+ \frac{\phi(\xi,\eta)}{D(\xi,\eta)} \right]^2 m_\alpha(\xi,\eta)\,.
 \end{split}
\end{align*}
Thus, using \eqref{malphasup} (with $k=2,0$), we see that
$|I|\leq\frac{C}{1-\alpha}|\eta|^{-2}$. 
An elementary but long computation shows that 
\begin{multline}
I\!I =\frac{m_\alpha(\xi,\eta)}{
 \eta^2  D^2(\xi,\eta)}\bigg[ 2\cosh\xi\frac{\sinh\xi}{\xi}
+ \cosh^2\xi-\frac{\sinh^2\xi}{\xi^2}\\
-\cosh^2\xi \frac{\sinh^2\eta}{\eta^2}
-2\cosh\xi \cosh\eta  \frac{\sinh\xi}{\xi}\frac{\sinh\eta}{\eta}   \bigg] \,. \label{++}
\end{multline}
which easily implies that
$|I\!I|\leq\frac{C}{1-\alpha}|\eta|^{-2}$. 
Therefore, 
$$
\big|\p^2_\eta m_{\alpha}(\xi,\eta)\big|
\le \frac{C}{1-\alpha}\frac{1}{|\eta|^2}
\,,
$$
for a constant $C$ independent of $\alpha$.

Finally, from \eqref{Derivative-eta}, with an elementary computation
we obtain that
\begin{align*}
  \p_\xi\p_\eta m_\alpha(\xi,\eta)
&=\p_\xi m_\alpha(\xi,\eta)\left[(\alpha-1) +\frac{\phi(\xi,\eta)}{D(\xi,\eta)}\right]\\
  &\qquad\qquad\qquad +\frac{m_\alpha(\xi,\eta)}{D^2(\xi,\eta)}
\bigg[\p_\xi \phi (\xi,\eta)D(\xi,\eta)-\phi(\xi,\eta) \p_\xi D(\xi,\eta)\bigg]\\
&=\p_\xi m_\alpha(\xi,\eta)
\left[(\alpha-1)+\frac{\phi(\xi,\eta)}{D(\xi,\eta)}\right]\\
&\qquad\qquad\qquad\qquad
+\frac{m_\alpha(\xi,\eta)}{D^2(\xi,\eta)}\bigg[\frac{1}{\xi\eta}
\Big(1-\cosh\xi\frac{\sinh\xi}{\xi}\Big)
\Big(1-\cosh\eta\frac{\sinh\eta}{\eta}\Big)\bigg].    
\end{align*}
Using \eqref{later-ref} and \eqref{der-xi} we now deduce  that
$$
|\p_\xi\p_\eta m_\alpha(\xi,\eta)|\leq
\frac{C}{1-\alpha}\, \frac{1}{|\xi\eta|}\,,
$$ 
for a constant $C$ independent of $\alpha$, and the lemma is proved. 
\qed
\ms

\subsection{Proof of Lemma \ref{Mar-lem-off-diag} }
  If $m_\alpha(\xi,\eta)$ denotes the function \eqref{malpha}, for all
  non-negative integer $k$ we
  clearly have
$$
|m_{\alpha,\gamma}(\xi,\eta)|\leq |m_\alpha(\xi,\eta)|
\leq\frac{C_k}{(1-\alpha)^{k+1}}\frac{1}{|\eta|^{k}} \,.
 $$
 By symmetry we also have the estimate
\begin{equation*}
|m_{\alpha,\gamma}(\xi,\eta)|\leq \frac{C_k}{(1-\gamma)^{k+1}}\frac{1}{|\xi|^k}. 
\end{equation*}
Therefore,
\begin{equation}\label{min-est-0}
 |m_{\alpha,\gamma}(\xi,\eta)|
\le C \min\Big(\frac{1}{1-\gamma},\frac{1}{1-\alpha}\Big)
\,.
 \end{equation}

Recalling that  $D(\xi,\eta)$ is given by
\eqref{Denominator}, arguing as in \eqref{Derivative-eta} we obtain
\begin{equation}\label{---}
\p_\eta
m_{\alpha,\gamma}(\xi,\eta)=\left[(\alpha-1)+\frac{\phi(\xi,\eta)}{D(\xi,\eta)}\right]m_{\alpha,\gamma}(\xi,\eta).
\end{equation}
Since $|m_{\alpha,\gamma}(\xi,\eta)|\leq |m_\alpha(\xi,\eta)|$, the
previous lemma guarantees that $|\p_\eta
m_{\alpha,\gamma}(\xi,\eta)|\leq \frac{C}{1-\alpha}|\eta|^{-1}$. 
Moreover,
\begin{align*}
|m_{\alpha,\gamma}(\xi,\eta)|
& \le \frac{e^{\gamma\xi}e^{\alpha\eta}}{\frac{\sinh\xi}{\xi} \cosh\eta}
  \le \frac{C}{(1-\gamma)|\eta|}
  \Big| \frac{ \eta e^{\alpha\eta}}{\cosh\eta}\Big| \\
& \le \frac{C}{(1-\gamma)(1-\alpha)}\frac{1}{|\eta|} \,.
\end{align*}
Inserting this and \eqref{min-est-0} into \eqref{---} we obtain 
\begin{align*}
|\p_\eta m_{\alpha,\gamma}(\xi,\eta)|
& \le  \frac{C}{(1-\gamma)}\frac{1}{|\eta|} +
\frac{C}{(1-\gamma)}\Big|\frac{\phi(\xi,\eta)}{D(\xi,\eta)}\Big| \\
& \le  \frac{C}{(1-\gamma)}\frac{1}{|\eta|} \,.
\end{align*}

Therefore,
$$
|\p_\eta
m_{\alpha,\gamma}(\xi,\eta)|\le 
C \min\Big(\frac{1}{1-\gamma},\frac{1}{1-\alpha}\Big) \frac{1}{|\eta|} \,.
$$

For the second order derivative with respect to $\eta$, similarly to
\eqref{2-Derivative-eta}, we obtain 
\begin{align*}
 \begin{split}
  \p^2_\eta & m_{\alpha,\gamma}(\xi,\eta)\\
  &=\left[(\alpha-1)+\frac{\phi(\xi,\eta)}{D(\xi,\eta)}\right]^2m_{\alpha,\gamma}(\xi,\eta)+\frac{m_{\alpha,\gamma}(\xi,\eta)}{D^2(\xi,\eta)}
\big[ \p_\eta \phi(\xi,\eta) D(\xi,\eta)
-\phi(\xi,\eta)\p_\eta D(\xi,\eta) \big] \,,
 \end{split}
\end{align*}
where $\phi$ is defined in \eqref{Derivative-eta}. 
Thus, using again the fact that $m_{\alpha,\gamma}$ is bounded from
above by $m_\alpha$, using the previous lemma,  we can easily conclude
(see \eqref{++})
that $|\p^2_\eta m_{\alpha,\gamma}(\xi,\eta)|\leq
\frac{C}{1-\alpha}|\eta|^{-2}$, but also 
$|\p^2_\eta m_{\alpha,\gamma}(\xi,\eta)|\leq
\frac{C}{1-\gamma}|\eta|^{-2}$; hence
$$
|\p^2_\eta
m_{\alpha,\gamma}(\xi,\eta)|\le 
C \min\Big(\frac{1}{1-\gamma},\frac{1}{1-\alpha}\Big) \frac{1}{|\eta|^2} \,.
$$

By symmetry we have the estimates 
\begin{equation*}
 |\p_\xi m_{\alpha,\gamma}(\xi,\eta)|
\leq \min\Big(\frac{1}{1-\gamma},\frac{1}{1-\alpha}\Big)  \frac{1}{|\xi|} \,,
\end{equation*}
and
\begin{equation*}
 |\p^2_\xi m_{\alpha,\gamma}(\xi,\eta)|
\leq\min\Big(\frac{1}{1-\gamma},\frac{1}{1-\alpha}\Big) \frac{1}{|\xi|^2}.
\end{equation*}
Finally, for the mixed derivative, we 
obtain
\begin{align*}
 \begin{split}
\p_\xi\p_\eta m_{\alpha,\gamma}(\xi,\eta)
&=\p_\xi m_{\alpha,\gamma}(\xi,\eta)\left[(\alpha-1) 
+\frac{\phi(\xi,\eta)}{D(\xi,\eta)}\right]\\
&\qquad\qquad\qquad
+\frac{m_{\alpha,\gamma}(\xi,\eta)}{D^2(\xi,\eta)} 
\bigg[\frac{1}{\xi\eta}\Big(1-\cosh\xi\frac{\sinh\xi}{\xi}\Big)
\Big(1-\cosh\eta\frac{\sinh\eta}{\eta}\Big)\bigg].
 \end{split}
 \end{align*}
The estimate 
\begin{equation*}
 |\p_\xi\p_\eta m_{\alpha,\gamma}(\xi,\eta)|
\leq C\min\bigg(\frac{1}{1-\gamma},\frac{1}{1-\alpha}\bigg) \frac{1}{|\xi\eta|}
\end{equation*}
now follows easily. 

Finally, we can conclude that $m_{\alpha,\gamma}$ is a Marcinkiewicz
multiplier whose norm of the associated multiplier operator is bounded
by a constant times $\min\big(\frac{1}{1-\gamma},\frac{1}{1-\alpha}\big)$. 
We conclude the proof by noticing that
$\min\big(\frac{1}{1-\gamma},\frac{1}{1-\alpha}\big)
\leq \frac{2}{(1-\alpha)+(1-\gamma)}$. 
\qed

\ms

\bibliographystyle{abbrv}
\bibliography{bibWormTopological}

\end{document}